\documentclass[thmsa,12pt]{article}
\usepackage{amssymb}

\usepackage{sw20lart}

\input{tcilatex}
\begin{document}

\begin{center}
\[
\begin{array}{cc}
\hspace{5cm}\hspace{1cm} & \text{\textit{In Memory of}} \\ 
\hspace{5cm}\hspace{1cm} & \text{\textit{Academician Marius Iosifescu}} \\ 
\hspace{5cm}\hspace{1cm} & \text{(1936}-\text{2020)}
\end{array}
\]

$G^{+}$ \textbf{Method in\ Action:}

\medskip

\textbf{New\ Classes\ of\ Nonnegative\ Matrices}

\medskip

\textbf{with Results}

\bigskip

\textbf{Udrea\ P\u {a}un}

\bigskip
\end{center}

\noindent The $G^{+}$ method is a new method, a powerful one, for the study
of (homogeneous and nonhomogeneous) products of nonnegative matrices --- for
problems on the products of nonnegative matrices. To study such products,
new classes of matrices are introduced: that of the sum-positive matrices,
that of the $\left[ \Delta \right] $-positive matrices on partitions (of the
column index sets), that of the $g_{k}^{+}$-matrices... On the other hand,
the $g_{k}^{+}$-matrices lead to necessary and sufficient conditions for the 
$k$-connected graphs. Using the $G^{+}$ method, we prove old and new results
(Wielandt Theorem and a generalization of it, etc.) on the products of
nonnegative matrices --- mainly, sum-positive, $\left[ \Delta \right] $%
-positive on partitions, irreducible, primitive, reducible, fully
indecomposable, scrambling, or Sarymsakov matrices, in some cases the
matrices being, moreover, $g_{k}^{+}$-matrices (not only irreducible).

\medskip


\medskip


\begin{center}
\medskip

\textbf{1. A SHORT INTRODUCTION}

\medskip
\end{center}

The $G$ method [17] led to good or very good results, see, \textit{e.g.},
[17, Theorems 1.6 and 1.8; see also Theorems 2.2 and 2.3], [18], [19], and
[20] --- see also [16; see, \textit{e.g.}, Examples 2.11, 2.19, and 2.22].

The $G^{+}$ method was suggested by the $G$ method and Theorem 2.2, both
from [17]. This method leads to good or very good results too; using this
method, we prove old and new results (Wielandt Theorem and a generalization
of it, etc.) on the (homogeneous or nonhomogeneous) products of nonnegative
matrices --- mainly, sum-positive, $\left[ \Delta \right] $-positive on
partitions (of the column index sets), irreducible, primitive, reducible,
fully indecomposable, scrambling, or Sarymsakov matrices, in some cases the
matrices being, moreover, $g_{k}^{+}$-matrices (not only irreducible).

The products of nonnegative matrices arise in several fields, such as,

1) (homogeneous and nonhomogeneous) Markov chains, see, \textit{e.g.}, [10]
and [21],

2) (fast or not too fast) exact sampling and other exactly solvable
problems/things (the computation of normalization constants, etc.), all
based on Markov chains (the methods used are not Monte Carlo), see, \textit{%
e.g.}, [18] and [19], where products of reducible stochastic matrices are
used (a surprise?),

3) Markov chain Monte Carlo, see, \textit{e.g.}, [5], [11], [12], and [13,
Chapter 9] --- since the Markov chain Monte Carlo methods are not too good,
it is very important to remove ``Monte Carlo'' from ``Markov chain Monte
Carlo'' in as many as possible cases, \textit{i.e.}, to obtain exact or
approximate good/efficient methods based on Markov chains, methods which are
not Monte Carlo (for exact methods, see 2) again), in as many as possible
cases,

4) probabilistic automata, see, \textit{e.g.}, [15],

5) fractals, see, \textit{e.g.}, [7, Sections 4.3 and 4.5] and [8, Chapter
11],

6) economics and some related fields, see, \textit{e.g.}, [8, Chapter 13]
and [9, pp. 487$-$488],

7) consensus, see, \textit{e.g.}, [21, pp. 153$-$158] and [23].

\begin{center}
\medskip

\textbf{2. }$G^{+}$ \textbf{METHOD --- DEFINITION AND\ BASIC\ RESULTS}

\medskip
\end{center}

In this section, we define the $G^{+}$ method and give certain basic results
of it. To define the $G^{+}$ method, we first define the sum-positive
matrices and $\left[ \Delta \right] $-positive matrices on partitions ---
the $\Delta $-positive matrices on partitions are also defined, and are also
important.

\smallskip

Set 
\[
\left\langle m\right\rangle =\left\{ 1,2,...,m\right\} \text{ (}m\in \Bbb{N},%
\text{ }m\geq 1\text{)}, 
\]
\[
\left\langle \left\langle m\right\rangle \right\rangle =\left\{
0,1,...,m\right\} \text{ (}m\in \Bbb{N}\text{)}, 
\]

\[
N_{m,n}=\left\{ P\left| \text{ }P\text{ is a nonnegative }m\times n\text{
matrix}\right. \right\} \text{,} 
\]
\[
S_{m,n}=\left\{ P\left| \text{ }P\text{ is a stochastic }m\times n\text{
matrix}\right. \right\} 
\]

\noindent --- here, a stochastic matrix is a row stochastic matrix ---,

\[
N_{n}=N_{n,n}\text{,} 
\]

\[
S_{n}=S_{n,n}. 
\]

The entry $\left( i,j\right) $ of a matrix $Z$ will be denoted $Z_{ij}$ or,
if confusion can arise, $Z_{i\rightarrow j}.$

\smallskip

Let $P=\left( P_{ij}\right) \in N_{m,n}$ ($i\in \left\langle m\right\rangle $
and $j\in \left\langle n\right\rangle $). $\left\langle m\right\rangle $ and 
$\left\langle n\right\rangle $ are called the \textit{index sets of} $P;$
moreover, $\left\langle m\right\rangle $ is called the \textit{row index set
of} $P$ while $\left\langle n\right\rangle $ is called the \textit{column
index set of} $P.$ If $P\in N_{n},$ $\left\langle n\right\rangle $ is called
the \textit{index set of} $P.$

\smallskip

Let $P=\left( P_{ij}\right) \in N_{m,n}$. Let $\emptyset \neq U\subseteq
\left\langle m\right\rangle $ and $\emptyset \neq V\subseteq \left\langle
n\right\rangle $. Set the matrices (these are submatrices of $P$) 
\[
P_{U}=\left( P_{ij}\right) _{i\in U,j\in \left\langle n\right\rangle },\text{
}P^{V}=\left( P_{ij}\right) _{i\in \left\langle m\right\rangle ,j\in V},%
\text{ and }P_{U}^{V}=\left( P_{ij}\right) _{i\in U,j\in V}. 
\]

\smallskip

\textit{Definition} 2.1. (See, \textit{e.g.}, [21, p. 80].) Let $P\in
N_{m,n}.$ We say that $P$ is a \textit{row-allowable matrix} if it has at
least one positive entry in each row. We say that $P$ is a \textit{%
column-allowable matrix} if it has at least one positive entry in each
column (equivalently, if the transpose of $P$ is row-allowable).

\smallskip

Below we define a central notion of this article.

\smallskip

\textit{Definition} 2.2. Let $P\in N_{m,n}$. Let $\emptyset \neq U\subseteq
\left\langle m\right\rangle $ and $\emptyset \neq V\subseteq \left\langle
n\right\rangle $. We say that $P$ is \textit{sum-positive on }$U\times V$ if 
\[
\sum\limits_{j\in V}P_{ij}>0,\text{ }\forall i\in U. 
\]

\smallskip

Based on Definition 2.2, we will use the generic name ``sum-positive
matrix/matrices''.

\smallskip

\textit{Remark} 2.3. $P$ is sum-positive on $U\times V$ if and only if $%
P_{U}^{V}$ is a row-allowable matrix, where $P\in N_{m,n},$ $\emptyset \neq
U\subseteq \left\langle m\right\rangle ,$ and $\emptyset \neq V\subseteq
\left\langle n\right\rangle $.

\smallskip

\textit{Remark} 2.4. Let $P\in N_{m,n}$.

(a) The column $j$ of $P$ is positive, \textit{i.e.}, $P^{\left\{ j\right\}
}>0\Longleftrightarrow $ $P$ is sum-positive on $\left\langle m\right\rangle
\times \left\{ j\right\} $.

(b) $P$ is positive $\Longleftrightarrow P^{\left\{ j\right\} }>0,$ $\forall
j\in \left\langle n\right\rangle \Longleftrightarrow $ $P$ is sum-positive
on $\left\langle m\right\rangle \times \left\{ j\right\} ,$ $\forall j\in
\left\langle n\right\rangle $.

\smallskip

\textit{Remark} 2.5. Let $P\in N_{m,n}.$

(a) If $P$ is positive, then it is sum-positive on $\left\langle
m\right\rangle \times \left\langle n\right\rangle $ --- more generally, on $%
U\times V,$ $\forall U,V,$ $\emptyset \neq U\subseteq \left\langle
m\right\rangle $ and $\emptyset \neq V\subseteq \left\langle n\right\rangle $%
. If $P$ is sum-positive on $\left\langle m\right\rangle \times \left\langle
n\right\rangle ,$ it does not follow that it is positive when $n\geq 2.$ Due
to these facts, the sum-positive matrices are generalizations of the
positive matrices.

(b) $P\neq 0\Longleftrightarrow \exists U\in \left\langle m\right\rangle ,$ $%
\exists V\in \left\langle n\right\rangle $ such that $P$ is sum-positive on $%
U\times V.$ Due to this fact, the class of all sum-positive matrices is
equal to the class of all nonnegative matrices without the zero matrices.

\smallskip

Let $\emptyset \neq U\subseteq \left\langle m\right\rangle $ and $\emptyset
\neq V\subseteq \left\langle n\right\rangle $. Set 
\[
G_{U,V}^{+}=G_{U,V}^{+}\left( m,n\right) =\left\{ P\left| \text{ }P\in
S_{m,n}\text{ and }P\text{ is sum-positive on }U\times V\right. \right\} 
\]

\noindent and 
\[
\overline{G}_{U,V}^{+}=\overline{G}_{U,V}^{+}\left( m,n\right) =\left\{
P\left| \text{ }P\in N_{m,n}\text{ and }P\text{ is sum-positive on }U\times
V\right. \right\} 
\]

\noindent --- obviously, 
\[
G_{U,V}^{+}\subseteq S_{m,n},\text{ }\overline{G}_{U,V}^{+}\subseteq N_{m,n},%
\text{ and }G_{U,V}^{+}\subseteq \overline{G}_{U,V}^{+}. 
\]

\smallskip

Below we give a basic result on the sum-positivity of matrices.

\smallskip

\textbf{THEOREM 2.6.} (i) \textit{Let} $P_{1}\in \overline{G}%
_{U_{1},U_{2}}^{+}\subseteq N_{m_{1},m_{2}}$ \textit{and }$P_{2}\in 
\overline{G}_{U_{2},U_{3}}^{+}\subseteq N_{m_{2},m_{3}}.$\textit{\ Then} 
\[
P_{1}P_{2}\in \overline{G}_{U_{1},U_{3}}^{+}\subseteq N_{m_{1},m_{3}}. 
\]

\smallskip

(ii) (\textit{a generalization of }(i)) \textit{Let }$P_{1}\in \overline{G}%
_{U_{1},U_{2}}^{+}\subseteq N_{m_{1},m_{2}},$\textit{\ }$P_{2}\in \overline{G%
}_{U_{2},U_{3}}^{+}\subseteq N_{m_{2},m_{3}},$ $...,$ $P_{n}\in \overline{G}%
_{U_{n},U_{n+1}}^{+}\subseteq N_{m_{n},m_{n+1}}.$\textit{\ Then} 
\[
P_{1}P_{2}...P_{n}\in \overline{G}_{U_{1},U_{n+1}}^{+}\subseteq
N_{m_{1},m_{n+1}}. 
\]

\smallskip

\textit{Proof}. (i) We have 
\[
\left( P_{1}P_{2}\right) _{U_{1}}^{U_{3}}=\left( P_{1}\right)
_{U_{1}}^{\left\langle m_{2}\right\rangle }\left( P_{2}\right)
_{\left\langle m_{2}\right\rangle }^{U_{3}}\geq \left( P_{1}\right)
_{U_{1}}^{U_{2}}\left( P_{2}\right) _{U_{2}}^{U_{3}}. 
\]

\noindent The matrices $\left( P_{1}\right) _{U_{1}}^{U_{2}}$ and $\left(
P_{2}\right) _{U_{2}}^{U_{3}}$ are row-allowable, so, $\left( P_{1}\right)
_{U_{1}}^{U_{2}}\left( P_{2}\right) _{U_{2}}^{U_{3}}$ is row-allowable, and,
further, $\left( P_{1}P_{2}\right) _{U_{1}}^{U_{3}}$ is row-allowable. By
Remark 2.3 we have 
\[
P_{1}P_{2}\in \overline{G}_{U_{1},U_{3}}^{+}. 
\]

(ii) Induction. $\square $

\smallskip

\textit{Remark} 2.7. Since $S_{m,n}\subseteq N_{m,n},$ $G_{U,V}^{+}\subseteq 
\overline{G}_{U,V}^{+},$ and a product of stochastic matrices is a
stochastic matrix, Theorem 2.6 holds, in particular, for stochastic
matrices. From, \textit{e.g.}, Theorem 2.6(i), we obtain

\textit{If} $P_{1}\in G_{U_{1},U_{2}}^{+}\subseteq S_{m_{1},m_{2}}$ \textit{%
and }$P_{2}\in G_{U_{2},U_{3}}^{+}\subseteq S_{m_{2},m_{3}},$\textit{\ then} 
\[
P_{1}P_{2}\in G_{U_{1},U_{3}}^{+}\subseteq S_{m_{1},m_{3}}. 
\]

\noindent For other results from this article, we can proceed similarly ---
we will omit to specify this fact further.

\smallskip

The next result is simple, beautiful, and important --- using it, it is
proved, \textit{e.g.}, the celebrated theorem of Wielandt for the index of
primitivity of a primitive matrix (for this theorem, see, \textit{e.g.}, [9,
p. 520] or, here, Theorem 4.24).

\smallskip

\textbf{THEOREM 2.8.} \textit{Let }$P_{1}\in \overline{G}_{U_{1},U_{2}}^{+}%
\subseteq N_{m_{1},m_{2}},$\textit{\ }$P_{2}\in \overline{G}%
_{U_{2},U_{3}}^{+}\subseteq N_{m_{2},m_{3}},$ $...,$ $P_{n}\in \overline{G}%
_{U_{n},U_{n+1}}^{+}\subseteq N_{m_{n},m_{n+1}}.$\textit{\ Let }$j\in
\left\langle m_{n+1}\right\rangle .$\textit{\ Suppose that }$%
U_{1}=\left\langle m_{1}\right\rangle $\textit{\ and }$U_{n+1}=\left\{
j\right\} .$\textit{\ Then} 
\[
\left( P_{1}P_{2}...P_{n}\right) ^{\left\{ j\right\} }>0 
\]

\noindent (\textit{i.e., the column }$j$\textit{\ of }$P_{1}P_{2}...P_{n}$%
\textit{\ is positive}).

\smallskip

\textit{Proof}. Remark 2.4(a) and Theorem 2.6(ii). $\square $

\smallskip

To give an example for the above result, we consider 
\[
P_{1}=\left( 
\begin{array}{cccc}
1 & 0 & 0 & 0 \\ 
0 & 2 & 0 & 0 \\ 
1 & 0 & 0 & 0 \\ 
0 & 1 & 0 & 3
\end{array}
\right) \text{ and }P_{2}=\left( 
\begin{array}{cccc}
0 & 0 & 1 & 0 \\ 
1 & 0 & 1 & 0 \\ 
0 & 0 & 2 & 0 \\ 
4 & 0 & 0 & 0
\end{array}
\right) . 
\]

\noindent Since $P_{1}\in \overline{G}_{\left\langle 4\right\rangle
,\left\langle 3\right\rangle }^{+}\subseteq N_{4}$ ($\left( P_{1}\right)
_{\left\langle 4\right\rangle }^{\left\langle 3\right\rangle }$ is
row-allowable) and $P_{2}\in \overline{G}_{\left\langle 3\right\rangle
,\left\{ 3\right\} }^{+}\subseteq N_{4}$, we have $\left( P_{1}P_{2}\right)
^{\left\{ 3\right\} }>0$ (\textit{i.e.}, the column 3 of $P_{1}P_{2}$ is
positive) --- this thing can also be obtained by direct computation.

\smallskip

Let $P\in N_{m,n}.$ Let $\emptyset \neq U\subseteq \left\langle
m\right\rangle $ and $\emptyset \neq V\subseteq \left\langle n\right\rangle $%
. Set 
\[
U\rightarrow V\text{ or, if confusion can arise, }U\stackrel{P}{\rightarrow }%
V 
\]

\noindent if $\forall i\in U,$ $\exists j\in V$ such that $P_{ij}>0.$ Set 
\[
V\leftarrow U\text{ or, if confusion can arise, }V\stackrel{P}{\leftarrow }U 
\]

\noindent if $U\rightarrow V$. Obviously, if $U\rightarrow V$ (or if $%
V\leftarrow U$), then $P$ is sum-positive on $U\times V.$

\smallskip

\textit{Remark} 2.9. Let $P_{1}\in N_{m_{1},m_{2}},$ $P_{2}\in
N_{m_{2},m_{3}},$ $...,$ $P_{n}\in N_{m_{n},m_{n+1}}.$\ Let $j\in
\left\langle m_{n+1}\right\rangle .$\ If 
\[
U_{1}=\left\langle m_{1}\right\rangle \rightarrow U_{2}\rightarrow
...\rightarrow U_{n}\rightarrow U_{n+1}=\left\{ j\right\} 
\]
\ \noindent (equivalently, if 
\[
U_{n+1}=\left\{ j\right\} \leftarrow U_{n}\leftarrow ...\leftarrow
U_{2}\leftarrow U_{1}=\left\langle m_{1}\right\rangle \text{),} 
\]
\noindent where, obviously, $\emptyset \neq U_{2}\subseteq \left\langle
m_{2}\right\rangle ,$ $...,$ $\emptyset \neq U_{n}\subseteq \left\langle
m_{n}\right\rangle ,$ then, by Theorem 2.8, 
\[
\left( P_{1}P_{2}...P_{n}\right) ^{\left\{ j\right\} }>0. 
\]

\noindent (We used 
\[
U_{1}=\left\langle m_{1}\right\rangle \rightarrow U_{2}\rightarrow
...\rightarrow U_{n}\rightarrow U_{n+1}=\left\{ j\right\} , 
\]

\noindent not 
\[
U_{1}=\left\langle m_{1}\right\rangle \stackrel{P_{1}}{\rightarrow }U_{2}%
\stackrel{P_{2}}{\rightarrow }...\stackrel{P_{n-1}}{\rightarrow }U_{n}%
\stackrel{P_{n}}{\rightarrow }U_{n+1}=\left\{ j\right\} , 
\]

\noindent because no confusion can arise.)

\smallskip

Set 
\[
\text{Par}\left( E\right) =\left\{ \Delta \left| \text{ }\Delta \text{ is a
partition of }E\right. \right\} , 
\]

\noindent where $E$ is a nonempty set. We will agree that the partitions do
not contain the empty set. $\left( E\right) $ is the improper (degenerate)
partition of $E.$

\smallskip

\textit{Definition}\textbf{\ }2.10. Let $\Delta _{1},\Delta _{2}\in $Par$%
\left( E\right) .$ We say that $\Delta _{1}$\textit{\ is finer than }$\Delta
_{2}$ if $\forall V\in \Delta _{1},$ $\exists W\in \Delta _{2}$ such that $%
V\subseteq W.$

\smallskip

Write $\Delta _{1}\preceq \Delta _{2}$ when $\Delta _{1}$ is finer than $%
\Delta _{2}.$

\smallskip

Set 
\[
\left( \left\{ i\right\} \right) _{i\in \left\{
s_{1},s_{2},...,s_{t}\right\} }=\left( \left\{ s_{1}\right\} ,\left\{
s_{2}\right\} ,...,\left\{ s_{t}\right\} \right) ; 
\]
\[
\left( \left\{ i\right\} \right) _{i\in \left\{
s_{1},s_{2},...,s_{t}\right\} }\in \text{Par}\left( \left\{
s_{1},s_{2},...,s_{t}\right\} \right) \text{ (}t\geq 1\text{).} 
\]

\noindent \textit{E.g.}, 
\[
\left( \left\{ i\right\} \right) _{i\in \left\langle n\right\rangle }=\left(
\left\{ 1\right\} ,\left\{ 2\right\} ,...,\left\{ n\right\} \right) . 
\]

\smallskip

Other generalizations of the positive matrices are given below, in
Definitions 2.11 and 2.12.

\smallskip

\textit{Definition} 2.11. Let $P\in N_{m,n}$. Let $\Delta \in $Par$\left(
\left\langle m\right\rangle \right) $ and $\Sigma \in $Par$\left(
\left\langle n\right\rangle \right) $. We say that $P$ is a $\left[ \Delta
\right] $\textit{-positive matrix on }$\Sigma $ if $\forall V\in \Sigma ,$ $%
\exists U\in \Delta $ such that $P$ is sum-positive on $U\times V$
(equivalently, if $\forall V\in \Sigma ,$ $\exists U\in \Delta $ such that $%
P_{U}^{V}$ is a row-allowable matrix --- see Remark 2.3). A $\left[ \Delta
\right] $-positive matrix on $\left( \left\{ j\right\} \right) _{j\in
\left\langle n\right\rangle }$ is called $\left[ \Delta \right] $\textit{%
-positive }for short.

\smallskip

\textit{Definition} 2.12. Let $P\in N_{m,n}$. Let $\Delta \in $Par$\left(
\left\langle m\right\rangle \right) $ and $\Sigma \in $Par$\left(
\left\langle n\right\rangle \right) $. We say that $P$ is a $\Delta $\textit{%
-positive matrix on }$\Sigma $ if $\Delta $ is the least fine partition for
which $P$ is a $\left[ \Delta \right] $-positive matrix on $\Sigma $. A $%
\Delta $-positive matrix on $\left( \left\{ j\right\} \right) _{j\in
\left\langle n\right\rangle }$ is called $\Delta $\textit{-positive} while a 
$\left( \left\langle m\right\rangle \right) $-positive matrix on $\Sigma $
is called \textit{positive\ on }$\Sigma $ for short. A positive matrix on $%
\left( \left\{ j\right\} \right) _{j\in \left\langle n\right\rangle }$ is
called \textit{positive} for short --- in this case, $P\in \overline{G}%
_{\left( \left\langle m\right\rangle \right) ,\left( \left\{ j\right\}
\right) _{j\in \left\langle n\right\rangle }}^{+},$ so, $P>0,$ and,
therefore, ``\textit{positive}'' is justified; we have more, namely, $P>0$
if and only if $P\in \overline{G}_{\left( \left\langle m\right\rangle
\right) ,\left( \left\{ j\right\} \right) _{j\in \left\langle n\right\rangle
}}^{+}$.

\smallskip

Based on Definition 2.11, we will use the generic name (warning!) ``$\left[
\Delta \right] $-positive matrix/matrices on partitions (of the column index
sets)'' --- this corresponds to the general case from Definition 2.11 ---
and ``$\left[ \Delta \right] $-positive matrix/matrices'' --- this
corresponds to the special case from Definition 2.11. For generic names
based on Definition 2.12, we proceed similarly.

\smallskip

The $\left[ \Delta \right] $-positive matrices on partitions are both
generalizations of the positive matrices --- see Definition 2.12 --- and of
the column-allowable ones --- see the next result (the collection of
column-allowable matrices and that of $\left[ \Delta \right] $-positive ones
are equal).

\smallskip

\textbf{THEOREM 2.13.} \textit{Let }$P\in N_{m,n}$ ($m,$ $n\geq 1$). \textit{%
Then }$P$\textit{\ is column-allowable if and only if }$\exists \Delta \in $%
Par$\left( \left\langle m\right\rangle \right) $ \textit{such that }$P$%
\textit{\ is }$\left[ \Delta \right] $\textit{-positive.}

\smallskip

\textit{Proof}. ``$\Longrightarrow $'' Definition 2.11, taking $\Delta
=\left( \left\{ i\right\} \right) _{i\in \left\langle m\right\rangle }.$

``$\Longleftarrow $'' Obvious (see Definition 2.11). $\square $

\smallskip

Let $\Delta \in $Par$\left( \left\langle m\right\rangle \right) $ and $%
\Sigma \in $Par$\left( \left\langle n\right\rangle \right) $. Set 
\[
G_{\Delta ,\Sigma }^{+}=G_{\Delta ,\Sigma }^{+}\left( m,n\right) =\left\{
P\left| \text{ }P\in S_{m,n}\text{ and }P\text{ is a }\left[ \Delta \right] 
\text{-positive matrix on }\Sigma \right. \right\} 
\]

\noindent and 
\[
\overline{G}_{\Delta ,\Sigma }^{+}=\overline{G}_{\Delta ,\Sigma }^{+}\left(
m,n\right) =\left\{ P\left| \text{ }P\in N_{m,n}\text{ and }P\text{ is a }%
\left[ \Delta \right] \text{-positive matrix on }\Sigma \right. \right\} 
\]

\noindent --- obviously, $G_{\Delta ,\Sigma }^{+}\subseteq S_{m,n},$ $%
\overline{G}_{\Delta ,\Sigma }^{+}\subseteq N_{m,n},$ and $G_{\Delta ,\Sigma
}^{+}\subseteq \overline{G}_{\Delta ,\Sigma }^{+}.$

\smallskip

When we study or even when we construct products of nonnegative matrices (in
particular, products of stochastic matrices) using $G_{U,V}^{+},$ $\overline{%
G}_{U,V}^{+},$ $G_{\Delta ,\Sigma }^{+},$ or $\overline{G}_{\Delta ,\Sigma
}^{+}$, we will refer this as the $G^{+}$\textit{\ method}. $G^{+}$ comes
from the verb \textit{to group} and its derivatives and the adjective 
\textit{positive}.

\smallskip

\textit{Remark} 2.14. The $G$ method from [17] can be renamed the $G^{s}$ 
\textit{method}, $s$ coming from the adjective \textit{stable}. But, for
simplification, we do not make this thing --- on the other hand, it is not
strictly necessary.

\smallskip

The next result is a basic one, and is somehow similar to Theorem 2.6.

\smallskip

\textbf{THEOREM 2.15.} (i) \textit{Let} $P_{1}\in \overline{G}_{\Delta
_{1},\Delta _{2}}^{+}\subseteq N_{m_{1},m_{2}}$ \textit{and }$P_{2}\in 
\overline{G}_{\Delta _{2},\Delta _{3}}^{+}\subseteq N_{m_{2},m_{3}}.$\textit{%
\ Then} 
\[
P_{1}P_{2}\in \overline{G}_{\Delta _{1},\Delta _{3}}^{+}\subseteq
N_{m_{1},m_{3}}. 
\]

\smallskip

(ii) (\textit{a generalization of }(i)) \textit{Let }$P_{1}\in \overline{G}%
_{\Delta _{1},\Delta _{2}}^{+}\subseteq N_{m_{1},m_{2}},$\textit{\ }$%
P_{2}\in \overline{G}_{\Delta _{2},\Delta _{3}}^{+}\subseteq
N_{m_{2},m_{3}}, $ $...,$ $P_{n}\in \overline{G}_{\Delta _{n},\Delta
_{n+1}}^{+}\subseteq N_{m_{n},m_{n+1}}.$\textit{\ Then} 
\[
P_{1}P_{2}...P_{n}\in \overline{G}_{\Delta _{1},\Delta _{n+1}}^{+}\subseteq
N_{m_{1},m_{n+1}}. 
\]

\smallskip

\textit{Proof}. (i) Let $W\in \Delta _{3}.$ Since $P_{2}\in \overline{G}%
_{\Delta _{2},\Delta _{3}}^{+},$ $\exists V\in \Delta _{2}$ such that $P_{2}$
is sum-positive on $V\times W$ --- equivalently, $\left( P_{2}\right)
_{V}^{W}$ is row-allowable. Since $P_{1}\in \overline{G}_{\Delta _{1},\Delta
_{2}}^{+},$ $\exists U\in \Delta _{1}$ such that $P_{1}$ is sum-positive on $%
U\times V$ --- equivalently, $\left( P_{1}\right) _{U}^{V}$ is
row-allowable. By Theorem 2.6(i), $P_{1}P_{2}$ is sum-positive on $U\times W$%
.

We conclude that 
\[
P_{1}P_{2}\in \overline{G}_{\Delta _{1},\Delta _{3}}^{+}. 
\]

(ii) Induction. $\square $

\smallskip

\textit{Remark} 2.16. Obviously, we have a similar remark to Remark 2.7 ---
this is left to the reader.

\smallskip

\smallskip

The next result is somehow similar to Theorem 2.8 and, from [17], to Theorem
1.6(i) --- it is also simple, beautiful, and important.

\smallskip

\textbf{THEOREM\ 2.17.} \textit{Let }$P_{1}\in \overline{G}_{\Delta
_{1},\Delta _{2}}^{+}\subseteq N_{m_{1},m_{2}},$\textit{\ }$P_{2}\in 
\overline{G}_{\Delta _{2},\Delta _{3}}^{+}\subseteq N_{m_{2},m_{3}},$ $...,$ 
$P_{n}\in \overline{G}_{\Delta _{n},\Delta _{n+1}}^{+}\subseteq
N_{m_{n},m_{n+1}}.$\textit{\ Suppose that }$\Delta _{1}=\left( \left\langle
m_{1}\right\rangle \right) $\textit{\ and }$\Delta _{n+1}=\left( \left\{
j\right\} \right) _{j\in \left\langle m_{n+1}\right\rangle }.$\textit{\ Then}
\[
P_{1}P_{2}...P_{n}>0. 
\]

\smallskip

\textit{Proof}. Definition 2.12 and Theorem 2.15(ii). $\square $

\bigskip

\begin{center}
\textbf{3. }$g_{k}^{+}$\textbf{-MATRICES AND }$k$\textbf{-CONNECTED\ GRAPHS}
\end{center}

\bigskip

In this section, we define the $g_{k}^{+}$-matrices. These matrices were
suggested by the application of $G^{+}$ method (when this method was
applied), and for them we give certain basic results and some examples.
Further, the $g_{k}^{+}$-matrices led to $k$-connected graphs --- necessary
and sufficient conditions are given for these graphs.

\smallskip

\textit{Definition} 3.1. Let $n\geq 2$ ($n\in \Bbb{N}$). Let $P\in N_{n}$.
Let $k\in \left\langle n-1\right\rangle $. We say that $P$ is a $g_{k}^{+}$%
\textit{-matrix} if $\forall F,$ $\emptyset \neq F\subset \left\langle
n\right\rangle ,$%
\[
\exists E,\text{ }\emptyset \neq E\subseteq F^{c},\text{ }\left| E\right|
\geq \min \left( k,\left| F^{c}\right| \right) ,\text{ and }P\in \overline{G}%
_{E,F}^{+} 
\]
($F^{c}$ is the complement of $F;$ $P\in \overline{G}_{E,F}^{+}$ (\textit{%
i.e.}, $P$ is sum-positive on $E\times F$) is equivalent to $P_{E}^{F}$ is a
row-allowable matrix (\textit{i.e.}, to $\forall i\in E,$ $\exists j\in F$
such that $P_{ij}>0$)).

\smallskip

Concerning the notion of $g_{k}^{+}$-matrix, $g$ and $+$ come from the $%
G^{+} $ method because the application of this method suggested the
consideration of $g_{k}^{+}$-matrices.

\smallskip

If $P\in N_{1}$ and $P\neq 0$ (equivalently, $P>0$), then, by definition, $P$
is a $g_{1}^{+}$-matrix.

\smallskip

Let $n\geq 2.$ Let $P\in N_{n}$. Let $\emptyset \neq F\subset \left\langle
n\right\rangle $. Set 
\[
D_{F}=D_{F}\left( P\right) =\left\{ i\left| \text{ }i\in F^{c}\text{ and }%
\exists j\in F\text{ such that }P_{ij}>0\right. \right\} . 
\]

\smallskip

\textbf{THEOREM 3.2. }\textit{Let }$n\geq 2$\textit{. Let }$P\in N_{n}$%
\textit{. Let }$k\in \left\langle n-1\right\rangle $\textit{. Then }$P$%
\textit{\ is a }$g_{k}^{+}$\textit{-matrix if and only if} 
\[
\left| D_{F}\right| \geq \min \left( k,\left| F^{c}\right| \right) ,\mathit{%
\ }\forall F,\mathit{\ }\emptyset \neq F\subset \left\langle n\right\rangle
. 
\]

\smallskip

\textit{Proof}. ``$\Longrightarrow $'' Let $\emptyset \neq F\subset
\left\langle n\right\rangle .$ By Definition 3.1, $\exists E,$ $\emptyset
\neq E\subseteq F^{c},$ $\left| E\right| \geq \min \left( k,\left|
F^{c}\right| \right) ,$ and $\forall i\in E,$ $\exists j\in F$ such that $%
P_{ij}>0.$ Obviously, $E\subseteq D_{F}$ ($D_{F}\subseteq F^{c}$); further,
we have $\left| D_{F}\right| \geq \left| E\right| $. So, $\left|
D_{F}\right| \geq \min \left( k,\left| F^{c}\right| \right) .$

\smallskip

``$\Longleftarrow $'' Obvious (taking $E=D_{F},$ $F$ fixed, $\emptyset \neq
F\subset \left\langle n\right\rangle $). $\square $

\smallskip

Set 
\[
G_{1,1}^{+}=\left\{ P\left| \text{ }P\in S_{1}\text{ and }P\text{ is a }%
g_{1}^{+}\text{-matrix}\right. \right\} , 
\]
\[
\overline{G}_{1,1}^{+}=\left\{ P\left| \text{ }P\in N_{1}\text{ and }P\text{
is a }g_{1}^{+}\text{-matrix}\right. \right\} , 
\]
\noindent and, for any $n\geq 2$ and $k\in \left\langle n-1\right\rangle ,$%
\[
G_{n,k}^{+}=\left\{ P\left| \text{ }P\in S_{n}\text{ and }P\text{ is a }%
g_{k}^{+}\text{-matrix}\right. \right\} 
\]

\noindent and 
\[
\overline{G}_{n,k}^{+}=\left\{ P\left| \text{ }P\in N_{n}\text{ and }P\text{
is a }g_{k}^{+}\text{-matrix}\right. \right\} . 
\]

\smallskip

\textit{Remark} 3.3. (a) Obviously, $G_{1,1}^{+}=\left\{ I\right\} ,$ $I$ is
the identity matrix, $I=\left( 1\right) ,$ $\overline{G}_{1,1}^{+}=\left\{
P\left| \text{ }P\in N_{1}\text{ and }P>0\right. \right\} ,$ and $%
G_{1,1}^{+}\subset \overline{G}_{1,1}^{+}.$

\smallskip

(b) Obviously, $G_{n,1}^{+}\supseteq G_{n,2}^{+}\supseteq ...\supseteq
G_{n,n-1}^{+}$ and $\overline{G}_{n,1}^{+}\supseteq \overline{G}%
_{n,2}^{+}\supseteq ...\supseteq \overline{G}_{n,n-1}^{+},$ $\forall n\geq 2$%
. Moreover, ``$\supseteq $'' can be replaced with ``$\supset $'' in all
places. For this, first, we consider the matrix $P\in N_{n}$ with $%
P^{\left\{ j\right\} }>0,$ $\forall j\in \left\langle n-1\right\rangle $ (%
\textit{i.e.}, the columns $1,2,...,n-1$ of $P$ are positive) and $%
P^{\left\{ n\right\} }$ has the first $k$ entries greater than $0$ and the
last $n-k$ entries equal to $0$ (therefore, $P_{1n},P_{2n},...,P_{kn}>0$ and 
$P_{k+1\rightarrow n}=P_{k+2\rightarrow n}=...=P_{n-1\rightarrow n}=P_{nn}=0$%
), where $k\in \left\langle n-2\right\rangle .$ Obviously, $P\in \overline{G}%
_{n,k}^{+}$ and $P\notin \overline{G}_{n,k+1}^{+},$ so, $\overline{G}%
_{n,k}^{+}\supset \overline{G}_{n,k+1}^{+}.$ For stochastic matrices, we can
proceed similarly.

\smallskip

(c) Obviously, $G_{n,1}^{+}\subset \overline{G}_{n,1}^{+},$ $%
G_{n,2}^{+}\subset \overline{G}_{n,2}^{+},$ $...,$ $G_{n,n-1}^{+}\subset 
\overline{G}_{n,n-1}^{+},$ $\forall n\geq 2.$

\smallskip

\textit{Definition} 3.4. (See, \textit{e.g.}, [9, p. 360].) Let $P\in N_{n}$%
. We say that $P$ is \textit{reducible} if either

(a) $n=1$ and $P=0$

\noindent or

(b) $n\geq 2,$ $\exists Q\in N_{n},$ $Q$ is a permutation matrix, and $%
\exists r\in \left\langle n-1\right\rangle $ such that 
\[
^{\text{t}}QPQ=\left( 
\begin{array}{cc}
X & 0 \\ 
Y & Z
\end{array}
\right) 
\]

\noindent ($^{\text{t}}Q$ is the transpose of $Q$), where $X\in N_{r},$ $0$
is a zero matrix, $0\in N_{r,n-r},$ $Y\in N_{n-r,r},$ and $Z\in N_{n-r,n-r}$
($X,$ $Y,$ and $Z$ can be zero matrices).

\smallskip

\textit{Definition} 3.5. (See, \textit{e.g.}, [9, p. 361].) Let $P\in N_{n}$%
. We say that $P$ is \textit{irreducible} if it is not reducible.

\smallskip

\textit{Remark} 3.6. Let $P\in N_{n}$. By Definitions 2.1 and 3.4, if $P$ is
not row-allowable, then it is reducible. By Definitions 2.1 and 3.5, if $P$
is irreducible, then it is row-allowable.

\smallskip

The next result is simple, but useful --- a basic result for the irreducible
matrices.

\smallskip

\textbf{THEOREM 3.7.} \textit{Let }$n\geq 2.$\textit{\ Let }$P\in N_{n}$%
\textit{. Suppose that }$P$\textit{\ is irreducible. Then the following
statements hold.}

(i) $\forall i\in \left\langle n\right\rangle ,$ $\exists j\in \left\langle
n\right\rangle ,$ $j\neq i,$\textit{\ such that }$P_{ij}>0.$

(ii) $\forall i\in \left\langle n\right\rangle ,$ $\exists k\in \left\langle
n\right\rangle ,$ $k\neq i,$\textit{\ such that }$P_{ki}>0.$

(iii) (\textit{a generalization of }(i)\textit{\ and} (ii)) $\forall A,$ $%
\emptyset \neq A\subset \left\langle n\right\rangle ,$ $\exists i\in A,$ $%
\exists j\in A^{c}$ \textit{such that} $P_{ij}>0.$

\smallskip

\textit{Proof}. (iii) Let $\emptyset \neq A\subset \left\langle
n\right\rangle .$ If $\forall i\in A,$ $\forall j\in A^{c},$ we have $%
P_{ij}=0,$ then, by Definition 3.4, $P$ is reducible. Contradiction. $%
\square $

\smallskip

\textit{Remark} 3.8. Let $n\geq 2.$ Let $P\in N_{n}$. If $\forall i\in
\left\langle n\right\rangle ,$ $\exists j,k\in \left\langle n\right\rangle ,$
$j\neq i,$ $k\neq i$ ($j=k$ or $j\neq k$) such that $P_{ij}>0$ and $P_{ki}>0$%
, it does not follow that $P$ is irreducible. Indeed, letting 
\[
P\in N_{4},\text{ }P=\left( 
\begin{array}{cccc}
1 & 1 & 0 & 0 \\ 
1 & 1 & 0 & 0 \\ 
1 & 0 & 0 & 1 \\ 
0 & 0 & 1 & 0
\end{array}
\right) , 
\]

\noindent $P$ has the above property, but it is reducible. Note that for $%
i=4 $, we have $j=k=3.$

\smallskip

In the next result we give a characterization of the irreducible matrices
using the $g_{1}^{+}$-matrices --- and thus we have examples of $g_{1}^{+}$%
-matrices.

\smallskip

\textbf{THEOREM 3.9. }\textit{Let }$P\in N_{n},$ $n\geq 1$\textit{. Then }$P$%
\textit{\ is irreducible if and only if }$P\in \overline{G}_{n,1}^{+}.$

\smallskip

\textit{Proof}. \textit{Case} 1. $n=1.$ No problem.

\textit{Case} 2. $n\geq 2.$

``$\Longrightarrow $'' Let $\emptyset \neq F\subset \left\langle
n\right\rangle .$ Since $P$ is irreducible and $\emptyset \neq F\subset
\left\langle n\right\rangle ,$ using Theorem 3.7(iii), $\exists i\in F^{c},$ 
$\exists j\in F$ such that $P_{ij}>0.$ We take $E=\left\{ i\right\} ,$ and
have $\left| E\right| =1$. So, $P$ is a $g_{1}^{+}$-matrix.

``$\Longleftarrow $'' Suppose that $P$ is not irreducible. It follows from
Definition 3.4 that $\exists Q\in N_{n},$ $Q$ is a permutation matrix, and $%
\exists r\in \left\langle n-1\right\rangle $ such that 
\[
^{\text{t}}QPQ=\left( 
\begin{array}{cc}
X & 0 \\ 
Y & Z
\end{array}
\right) , 
\]

\noindent where $X\in N_{r},$ $0$ is a zero matrix, $0\in N_{r,n-r},$ $Y\in
N_{n-r,r},$ and $Z\in N_{n-r,n-r}.$ Suppose that the rows of $X$ and of $0$
are $i_{1},i_{2},,...,i_{r}$ --- the columns of $X$ are $%
i_{1},i_{2},,...,i_{r}$ --- and that the columns of $0$ are $%
i_{r+1},i_{r+2},,...,i_{n}.$ Take $F=\left\{
i_{r+1},i_{r+2},,...,i_{n}\right\} .$ Then $\forall E,$ $\emptyset \neq
E\subseteq F^{c}=\left\{ i_{1},i_{2},,...,i_{r}\right\} ,$ $\forall i\in E,$ 
$\forall j\in F$ we have $P_{ij}=0.$ So, $P\notin \overline{G}_{n,1}^{+}.$
Contradiction. $\square $

\smallskip

\textbf{THEOREM 3.10. }\textit{Let }$P\in N_{n}$\textit{. If }$P\in 
\overline{G}_{n,k}^{+}$\textit{\ for some }$k\in \left\langle
n-1\right\rangle ,$\textit{\ then }$P$\textit{\ is row-allowable.}

\smallskip

\textit{Proof}. Definition 2.1, Remarks 3.3(b) and 3.6, and Theorem 3.9 (any
irreducible matrix is row-allowable). $\square $

\smallskip

\textit{Remark} 3.11. If a nonnegative $n\times n$ matrix is row-allowable,
it does not follow that it is a $g_{k}^{+}$-matrix for some $k\in
\left\langle n-1\right\rangle $ if $n\geq 2.$ Any nonnegative $n\times n$
matrix which is row-allowable and reducible is not a $g_{k}^{+}$-matrix, $%
\forall k\in \left\langle n-1\right\rangle $. \textit{E.g.}, let 
\[
P\in N_{n},\text{ }P=\left( 
\begin{array}{ccccc}
1 & 0 & 0 & \cdots & 0 \\ 
1 & 0 & 0 & \cdots & 0 \\ 
\vdots & \vdots & \vdots & \cdots & \vdots \\ 
1 & 0 & 0 & \cdots & 0
\end{array}
\right) ,\text{ }n\geq 2. 
\]

\noindent $P$ is row-allowable (it is reducible), but $P\notin \overline{G}%
_{n,k}^{+},$ $\forall k\in \left\langle n-1\right\rangle .$

\smallskip

\textbf{THEOREM 3.12.} \textit{Let }$n\geq 2$\textit{\ and }$k\in
\left\langle n-1\right\rangle .$\textit{\ Let }$P\in N_{n}$\textit{. If }$%
P\in \overline{G}_{n,k}^{+},$\textit{\ then the matrix }$P_{\left\langle
n\right\rangle -\left\{ j\right\} }^{\left\{ j\right\} }$\textit{\ has at
least }$k$\textit{\ positive entries, }$\forall j\in \left\langle
n\right\rangle $\textit{\ }(\textit{i.e., in each column, }$P$\textit{\ has
at least }$k$\textit{\ positive entries which are not situated on the main
diagonal})\textit{.}

\smallskip

\textit{Proof}. Let $j\in \left\langle n\right\rangle .$ Take $F=\left\{
j\right\} .$ Since $P\in \overline{G}_{n,k}^{+},$ by Definition 3.1, 
\[
\exists E,\text{ }\emptyset \neq E\subseteq \left\langle n\right\rangle
-\left\{ j\right\} ,\text{ }\left| E\right| \geq \min \left( k,\left|
\left\langle n\right\rangle -\left\{ j\right\} \right| \right) ,\text{ and }%
P_{ij}>0,\text{ }\forall i\in E. 
\]
\noindent We have 
\[
\min \left( k,\left| \left\langle n\right\rangle -\left\{ j\right\} \right|
\right) =k 
\]
\noindent because $n\geq 2,$ $k\in \left\langle n-1\right\rangle ,$ and $%
\left| \left\{ j\right\} \right| =1.$ It follows that $\left| E\right| \geq
k.$ Since $P_{ij}>0,$ $\forall i\in E,$ further, it follows that, in column $%
j,$ $P$\textit{\ }has at least $k$ positive entries which are not situated
on the main diagonal. $\square $

\smallskip

\textit{Remark} 3.13. Let $n\geq 2$ and $k\in \left\langle n-1\right\rangle
. $ If a nonnegative $n\times n$ matrix has, in each column, at least $k$
positive entries which are not situated on the main diagonal, it does not
follow that it is a $g_{k}^{+}$-matrix. \textit{E.g.}, let 
\[
P\in N_{6},\text{ }P=\left( 
\begin{array}{cccccc}
0 & * & * & 0 & 0 & * \\ 
\ast & 0 & * & 0 & 0 & 0 \\ 
\ast & * & 0 & 0 & 0 & 0 \\ 
\ast & 0 & 0 & 0 & * & * \\ 
0 & 0 & 0 & * & 0 & * \\ 
0 & 0 & 0 & * & * & 0
\end{array}
\right) , 
\]

\noindent where ``$*$'' stands for a positive entry. $P$ has, in each
column, at least 2 positive entries which are not situated on the main
diagonal. But $P\notin \overline{G}_{6,2}^{+}.$ Indeed, if we take $%
F=\left\{ 1,2,3\right\} $ --- see Definition 3.1 ---, $E$ can be ($\left|
E\right| \geq 2$) $\left\{ 4,5\right\} ,$ $\left\{ 4,6\right\} ,$ $\left\{
5,6\right\} ,$ or $\left\{ 4,5,6\right\} $, but $P\notin \overline{G}%
_{E,F}^{+},$ $\forall E,$ $E$ being $\left\{ 4,5\right\} ,$ $\left\{
4,6\right\} ,$ $\left\{ 5,6\right\} ,$ or $\left\{ 4,5,6\right\} .$

\smallskip

By Remark 3.3(b) and Theorem 3.9 we know that any $g_{k}^{+}$-matrix is
irreducible; any irreducible matrix is a $g_{1}^{+}$-matrix --- and thus we
have examples of $g_{1}^{+}$-matrices. Further, we give other examples of $%
g_{k}^{+}$-matrices for $k=1$ or, specially, for $k\geq 2.$

\smallskip

\textit{Example} 3.14. Let $n\geq 2.$ Let $P\in N_{n}$. If $P>0$ or, more
generally, if $P$ is a matrix with $P_{ij}>0,$ $\forall i,j\in \left\langle
n\right\rangle ,$ $i\neq j,$ then it is a $g_{k}^{+}$-matrix, $\forall k\in
\left\langle n-1\right\rangle .$ Moreover, $P$ is a matrix with $P_{ij}>0,$ $%
\forall i,j\in \left\langle n\right\rangle ,$ $i\neq j,$ if and only if it
is a $g_{n-1}^{+}$-matrix.

\smallskip

The irreducible nonnegative matrices are (either) aperiodic or periodic,
see, \textit{e.g.}, [10, pp. 52$-$53], [21, Sections 1.2 and 1.3 (p. 18,
etc.)], and, here, Definition 4.20.

\smallskip

\textit{Example} 3.15. Consider the periodic irreducible matrix with period $%
t$ ($t\geq 2$) 
\[
P=\left( 
\begin{array}{ccccc}
& Q_{1} &  &  &  \\ 
&  & Q_{2} &  &  \\ 
&  &  & \ddots &  \\ 
&  &  &  & Q_{t-1} \\ 
Q_{t} &  &  &  & 
\end{array}
\right) \in N_{n}, 
\]

\noindent where $Q_{1}\in N_{n_{1},n_{2}},$ $Q_{2}\in N_{n_{2},n_{3}},$ $%
..., $ $Q_{t}\in N_{n_{t},n_{1}}$ $n_{1},n_{2},...,n_{t}\geq 1$ ($%
n_{1}+n_{2}+...+n_{t}=n$). Suppose that $Q_{1},$ $Q_{2},$ $...,$ $Q_{t}>0.$
Let $1\leq k\leq m=\min\limits_{1\leq l\leq t}n_{l}$ (obviously, in this
case, $k\leq n-1$). Then $P\in \overline{G}_{n,k}^{+},$ $\forall k\in
\left\langle m\right\rangle $ --- it is easy to see this, it is also easy to
see that $P\notin \overline{G}_{n,m+1}^{+}$ (but $P\in \overline{G}%
_{n,m}^{+} $).

\smallskip

\textit{Definition} 3.16. (See, \textit{e.g.}, [9, p. 356] and [10, p.
222].) Let $P\in N_{m,n}$. Set 
\[
\overline{P}=\left( \overline{P}_{ij}\right) \in N_{m,n},\text{ }\overline{P}%
_{ij}=\left\{ 
\begin{array}{l}
1\text{ if }P_{ij}>0, \\ 
0\text{ if }P_{ij}=0,
\end{array}
\right. 
\]

\noindent $\forall i\in \left\langle m\right\rangle ,$ $\forall j\in
\left\langle n\right\rangle .$ We call $\overline{P}$ the \textit{indicator }%
(or \textit{incidence})\textit{\ matrix of }$P$.

\smallskip

\textit{Definition} 3.17. (See, \textit{e.g.}, [8, p. 61].) Let $P\in
N_{m,n} $. Let $B\in N_{m,n}$ be a $\left( 0,1\right) $-matrix, \textit{i.e.}%
, $B$ is a matrix with $B_{ij}\in \left\{ 0,1\right\} ,$ $\forall i\in
\left\langle m\right\rangle ,$ $\forall j\in \left\langle n\right\rangle .$
We say that $P$ \textit{has the pattern} $B$ if $\overline{P}\geq B.$

\smallskip

\textit{Example} 3.18. Consider an aperiodic irreducible matrix $Q\in N_{n}$
(in this case, $Q$ is a primitive matrix --- see, in Section 4, Theorem
4.21). Consider that $Q$ has the pattern $\overline{P},$ where $P$ is the
matrix from Example 3.15. \textit{E.g.}, $Q=P+A$ is an aperiodic irreducible
matrix and has the pattern $\overline{P}$ if, \textit{e.g.}, $A\in N_{n}$
and $A_{11}>0$ or $A=I\in N_{n},$ $I$ is identity matrix. Then --- see
Example 3.15 --- $Q\in \overline{G}_{n,k}^{+},$ $\forall k\in \left\langle
m\right\rangle ,$ and $Q\notin \overline{G}_{n,m+1}^{+}.$

\smallskip

Below we will characterize the $k$-connected graphs by means of the $%
g_{k}^{+}$-matrices --- necessary and sufficient conditions are given for
these graphs. To do this thing, we need certain notions from the graph
theory. For the graph theory (basic notions, notation, etc.), see, \textit{%
e.g.}, [2]-[4], [6], and [22] --- all or some of these references (all are
book references) could be available; Wikipedia could also be useful...

\smallskip

We work with nondirected finite graphs. Moreover, we work with graphs
without multiple edges, but the loops (not the multiple loops) are allowed.

\smallskip

\textit{Definition} 3.19. (See, \textit{e.g.}, [9, p. 168].) Let $\mathcal{G}%
=\left( \mathcal{V},\mathcal{E}\right) $ be a (nondirected finite) graph
(without multiple edges, the loops are allowed), where $\mathcal{V}=\left\{
V_{1},V_{2},...,V_{n}\right\} $ is the vertex set ($n\geq 1$) and $\mathcal{E%
}$ is the edge set ($\left| \mathcal{E}\right| \geq 0$). Set 
\[
A=\left( A_{ij}\right) \in N_{n},\text{ }A_{ij}=\left\{ 
\begin{array}{ll}
1\smallskip & \text{if }\left[ V_{i},V_{j}\right] \in \mathcal{E}, \\ 
0 & \text{if }\left[ V_{i},V_{j}\right] \notin \mathcal{E},
\end{array}
\right. 
\]

\noindent $\forall i,j\in \left\langle n\right\rangle .$ $A$ is called the 
\textit{adjacency matrix of }$\mathcal{G}$\textit{.}

\textit{\smallskip }

Let $\mathcal{G}=\left( \mathcal{V},\mathcal{E}\right) $ be a graph. Let $%
\mathcal{W}\subset \mathcal{V}$. Consider the graph $\mathcal{G}-\mathcal{W}$%
. If $\mathcal{W}\neq \emptyset ,$ $\mathcal{G}-\mathcal{W}$ is the
(sub)graph obtained from $\mathcal{G}$ by deleting the vertices in $\mathcal{%
W}$ together with, if any, their incident edges. $\mathcal{G}-\mathcal{W}=%
\mathcal{G}$ if $\mathcal{W}=\emptyset .$ (See, \textit{e.g.}, [3, p. 9].)

\smallskip

\textit{Definition} 3.20. (See, \textit{e.g.}, [3, p. 42].) Let $\mathcal{G}%
=\left( \mathcal{V},\mathcal{E}\right) $ be a graph with $n$ vertices, $%
n\geq 1.$ Let $\mathcal{C}\subset \mathcal{V}$ with $\left| \mathcal{C}%
\right| \geq 0.$ We say that $\mathcal{C}$ is a \textit{vertex cut of }$%
\mathcal{G}$ if $\mathcal{G}-\mathcal{C}$ is disconnected ($\mathcal{G}-%
\mathcal{C}$ is disconnected $\Longrightarrow $ $n\geq 2,$ $\left| \mathcal{V%
}-\mathcal{C}\right| \geq 2,$ and $\left| \mathcal{C}\right| \leq n-2$; when 
$n=1,$ the graph $\mathcal{G}$ has no vertex cuts ($\emptyset $ is not a
vertex cut); when $n\geq 2$ and, moreover, the graph $\mathcal{G}$ is
disconnected, $\emptyset $ is a vertex cut). We say that the vertex cut $%
\mathcal{C}$ is a $k$-\textit{vertex cut }if $\left| \mathcal{C}\right| =k$
(obviously, now is very obviously, $0\leq k\leq n-2$).

\smallskip

$\mathcal{K}_{n}$ is the complete graph with $n$ vertices, $n\geq 1.$ $%
\mathcal{K}_{1}$ has one vertex and no edge. $\cong $ is the isomorphism
relation for graphs and $\ncong $ is its negation.

\smallskip

\textbf{THEOREM 3.21.} \textit{Let }$\mathcal{G}=\left( \mathcal{V},\mathcal{%
E}\right) $\textit{\ be a graph with }$n$\textit{\ vertices, }$n\geq 1.$%
\textit{\ Let }$\mathcal{G}^{\prime }$\textit{\ be the graph obtained from }$%
\mathcal{G}$\textit{\ deleting, if any, the loops. Then }$\mathcal{G}%
^{\prime }\ncong \mathcal{K}_{n}$\textit{\ if and only if }$\exists \mathcal{%
C}$\textit{, }$\mathcal{C}$\textit{\ is a vertex cut of }$\mathcal{G}$%
\textit{.}

\smallskip

\textit{Proof}. \textit{Case} 1. $n=1.$ Nothing to prove.

\textit{Case} 2. $n\geq 2.$

``$\Longrightarrow $'' Since $\mathcal{G}^{\prime }\ncong \mathcal{K}_{n},$
it follows that $\exists V,W\in \mathcal{G},$ $V\neq W,$ such that $V$ and $%
W $ are not adjacent. Let 
\[
\mathcal{C}=\left\{ X\left| \text{ }X\in \mathcal{V}\text{ and }X\text{ is
adjacent to }V\right. \right\} . 
\]
\noindent We can have $\mathcal{C}=\emptyset $ or $\mathcal{C}\neq \emptyset 
$ --- \textit{e.g.}, $\mathcal{C}=\emptyset $ when $n=2.$ If $\mathcal{C}%
=\emptyset ,$ the graph $\mathcal{G}^{\prime }$ is disconnected while, if $%
\mathcal{C}\neq \emptyset ,$ the graph $\mathcal{G}^{\prime }$ is
disconnected or connected. So, both when $\mathcal{C}=\emptyset $ and when $%
\mathcal{C}\neq \emptyset ,$ the graph $\mathcal{G}^{\prime }-\mathcal{C}$
is disconnected. It follows that $\mathcal{G}-\mathcal{C}$ is disconnected.
So, $\mathcal{C}$ is a vertex cut of $\mathcal{G}$.

\smallskip

``$\Longleftarrow $'' Obvious (because $\mathcal{G}-\mathcal{C}$ is
disconnected, so, $\exists Y,Z\in \mathcal{G},$ $Y\neq Z,$ such that $Y$ and 
$Z$ are not adjacent). $\square $

\smallskip

Due to Theorem 3.21, the next definition makes sense.

\smallskip

\textit{Definition} 3.22. (See, \textit{e.g.}, [3, p. 42].) Let $\mathcal{G}%
=\left( \mathcal{V},\mathcal{E}\right) $ be a graph with $n$ vertices, $%
n\geq 1.$ Let $\mathcal{G}^{\prime }$ be the graph obtained from $\mathcal{G}
$ deleting, if any, the loops. Set 
\[
\mathtt{k}\left( \mathcal{G}\right) =\left\{ 
\begin{array}{l}
\text{the minimum }k\text{ for which }\mathcal{G}\text{ has a }k\text{%
-vertex cut if }\mathcal{G}^{\prime }\ncong \mathcal{K}_{n},\medskip \\ 
\left| \mathcal{V}\right| -1\text{ (equivalently, }n-1\text{) if }\mathcal{G}%
^{\prime }\cong \mathcal{K}_{n}.
\end{array}
\right. 
\]

\noindent $\mathtt{k}\left( \mathcal{G}\right) $ is called the \textit{%
connectivity of }$\mathcal{G}$.

\smallskip

\textit{Definition} 3.23. (See, \textit{e.g.}, [3, p. 42].) Let $\mathcal{G}%
=\left( \mathcal{V},\mathcal{E}\right) $ be a graph with $n$ vertices, $%
n\geq 1.$ Let $k\geq 0.$ We say that $\mathcal{G}$ is $k$\textit{-connected}
if $\mathtt{k}\left( \mathcal{G}\right) \geq k$ ($\mathtt{k}\left( \mathcal{G%
}\right) \geq k$ implies $k\leq n-1,$ so, $k\in \left\langle \left\langle
n-1\right\rangle \right\rangle $).

\smallskip

\textit{Remark} 3.24. (a) $0\leq \mathtt{k}\left( \mathcal{G}\right) \leq
n-1 $ for any graph $\mathcal{G}$ with $n$ vertices, $n\geq 1;$ $\mathtt{k}%
\left( \mathcal{K}_{1}\right) =0,$ $\mathtt{k}\left( \mathcal{K}_{n}\right)
=n-1.$

\smallskip

(b) By Definition 3.23 and (a) any graph (with $n$ vertices, $n\geq 1$) is $%
0 $-connected, and, conversely, any $0$-connected graph is a graph. Any
graph with one vertex or which is disconnected has the connectivity equal to 
$0,$ so, it is only $0$-connected.

\smallskip

(c) Any connected graph is $1$-connected, and, conversely, any $1$-connected
graph is a connected graph.

\noindent (For $0$-connected graphs, see, \textit{e.g.}, [6, p. 11] (for $k$%
-connected graphs, see Section 1.4 and Chapter 3).)

\smallskip

The next result is a bridge between the $k$-connected graphs and $g_{k}^{+}$%
-matrices --- and thus we have other examples of $g_{k}^{+}$-matrices.

\smallskip

\textbf{THEOREM 3.25.} \textit{Let }$\mathcal{G}=\left( \mathcal{V},\mathcal{%
E}\right) $\textit{\ be a graph with }$n$\textit{\ vertices, }$n\geq 2.$%
\textit{\ Let }$k\in \left\langle n-1\right\rangle $ (\textit{this implies} $%
n\geq k+1$).\textit{\ Then }$\mathcal{G}$\textit{\ is }$k$\textit{-connected
if and only if its adjacency matrix, }$A,$\textit{\ is a }$g_{k}^{+}$\textit{%
-matrix.}

\smallskip

\textit{Proof}. \textit{Case} 1. $k=1.$ Theorem 3.9, and the fact that $%
\mathcal{G}$ is connected/$1$-connected if and only if its adjacency matrix
is irreducible.

\smallskip

\textit{Case} 2. $k\neq 1.$ Suppose that $\mathcal{V}=\left\{
V_{1},V_{2},...,V_{n}\right\} .$ Taking into account Definition 3.19, we
consider the bijective function $f:\left\langle n\right\rangle
\longrightarrow \mathcal{V},$ $i\longmapsto f\left( i\right) =V_{i}.$ Let $%
\emptyset \neq T\subset \left\langle n\right\rangle .$ Set 
\[
\mathcal{V}_{T}=\left\{ V_{i}\left| \text{ }i\in T\text{ and }f\left(
i\right) =V_{i}\right. \right\} . 
\]

\noindent Obviously, $\emptyset \neq \mathcal{V}_{T}\subset \mathcal{V}.$

``$\Longrightarrow $'' Suppose that $A$ is not a $g_{k}^{+}$-matrix. By
Theorem 3.2, $\exists F,$ $\emptyset \neq F\subset \left\langle
n\right\rangle ,$ such that $\left| D_{F}\right| <\min \left( k,\left|
F^{c}\right| \right) $ (recall that $k\geq 2$). Since $\mathcal{G}$ is $k$%
-connected, it is $1$-connected, so, $A$ is irreducible. It follows that $%
\left| D_{F}\right| \geq 1.$ From $\left| D_{F}\right| <\min \left( k,\left|
F^{c}\right| \right) ,$ we have $\left| D_{F}\right| <\left| F^{c}\right| .$
Therefore, we have $\emptyset \neq D_{F}\subset F^{c}.$ By the definition of 
$D_{F}$ and the fact that $F^{c}-D_{F}\neq \emptyset ,$ we have $A_{ij}=0,$ $%
\forall i\in F^{c}-D_{F},$ $\forall j\in F$ --- equivalently, we have $%
\left[ V_{i},V_{j}\right] \notin \mathcal{E}$, $\forall i\in F^{c}-D_{F},$ $%
\forall j\in F.$ It follows that $\mathcal{V}_{D_{F}}$ is a vertex cut of $%
\mathcal{G}$ ($\mathcal{G}-\mathcal{V}_{D_{F}}$ is disconnected). Since $%
\mathcal{G}$ is $k$-connected, $\mathtt{k}\left( \mathcal{G}\right) \geq k.$
By $k>\left| D_{F}\right| $ and $\left| D_{F}\right| =\left| \mathcal{V}%
_{D_{F}}\right| ,$ we obtain $\mathtt{k}\left( \mathcal{G}\right) >\left| 
\mathcal{V}_{D_{F}}\right| .$ Contradiction. So, $A$ is a $g_{k}^{+}$-matrix.

\smallskip

``$\Longleftarrow $'' Suppose that $\mathcal{G}$ is not $k$-connected. It
follows that $\mathtt{k}\left( \mathcal{G}\right) <k.$ By Remark 3.3(b), $A$
is a $g_{1}^{+}$-matrix. By Theorem 3.9, $A$ is irreducible. It follows that 
$\mathcal{G}$ is connected/$1$-connected. So, $\mathtt{k}\left( \mathcal{G}%
\right) \geq 1.$ We now have $1\leq \mathtt{k}\left( \mathcal{G}\right) <k$ (%
$k\geq 2$). It follows that $\exists \mathcal{C},$ $\mathcal{C}$ is a vertex
cut of $\mathcal{G}$, with $1\leq \left| \mathcal{C}\right| <k.$ Further, it
follows that $\mathcal{G}-\mathcal{C}$ is disconnected, and, as a result, $%
\left| \mathcal{V}-\mathcal{C}\right| \geq 2$ --- further, $\exists \mathcal{%
W}_{1},\mathcal{W}_{2}\subset \mathcal{V}-\mathcal{C},$ $\mathcal{W}_{1},%
\mathcal{W}_{2}\neq \emptyset ,$ $\mathcal{W}_{1}\cap \mathcal{W}%
_{2}=\emptyset ,$ $\mathcal{W}_{1}\cup \mathcal{W}_{2}=\mathcal{V}-\mathcal{C%
},$ $1\leq \left| \mathcal{W}_{1}\right| \leq \left| \mathcal{W}_{2}\right| $%
, and $\left[ U,V\right] \notin \mathcal{E},$ $\forall U\in \mathcal{W}_{1},$
$\forall V\in \mathcal{W}_{2}.$ Let 
\[
F=\left\{ i\left| \text{ }V_{i}\in \mathcal{W}_{1}\right. \right\} . 
\]
\[
\text{(}F=\left\{ i\left| \text{ }V_{i}\in \mathcal{W}_{2}\right. \right\} 
\text{ is also good.)} 
\]

\noindent Obviously, $\emptyset \neq F\subset \left\langle n\right\rangle .$
Since $\mathcal{C}\neq \emptyset $ (because $\left| \mathcal{C}\right| \geq
1 $) and $\left[ U,V\right] \notin \mathcal{E},$ $\forall U\in \mathcal{W}%
_{1}, $ $\forall V\in \mathcal{W}_{2},$ we have 
\[
D_{F}\subseteq \left\{ i\left| \text{ }V_{i}\in \mathcal{C}\right. \right\}
. 
\]

\noindent So, $\left| D_{F}\right| \leq \left| \mathcal{C}\right| .$ Since $%
\left| \mathcal{C}\right| <k,$ we have $\left| D_{F}\right| <k.$ On the
other hand, we have $\left| D_{F}\right| <\left| F^{c}\right| $ because 
\[
F^{c}=\left\{ i\left| \text{ }V_{i}\in \mathcal{C}\cup \mathcal{W}%
_{2}\right. \right\} =\left\{ i\left| \text{ }V_{i}\in \mathcal{C}\right.
\right\} \cup \left\{ i\left| \text{ }V_{i}\in \mathcal{W}_{2}\right.
\right\} \supset 
\]
\[
\supset \left\{ i\left| \text{ }V_{i}\in \mathcal{C}\right. \right\}
\supseteq D_{F}. 
\]

\noindent From $\left| D_{F}\right| <k$ and $\left| D_{F}\right| <\left|
F^{c}\right| ,$ we have $\left| D_{F}\right| <\min \left( k,\left|
F^{c}\right| \right) .$ Contradiction (because $A$ is a $g_{k}^{+}$-matrix).
So, $\mathcal{G}$ is $k$-connected. $\square $

\smallskip

Let $\mathcal{G}=\left( \mathcal{V},\mathcal{E}\right) $ be a graph with $n$
vertices, $n\geq 2.$ Let $\emptyset \neq \mathcal{Y}\subset \mathcal{V}.$
Set 
\[
D_{\mathcal{Y}}=\left\{ V\left| V\in \mathcal{Y}^{c}\text{ and }\exists W\in 
\mathcal{Y}\text{ such that }\left[ V,W\right] \in \mathcal{E}\right.
\right\} . 
\]

\smallskip

We arrived at an interesting result on the $k$-connected graphs.

\smallskip

\textbf{THEOREM 3.26. }\textit{Let }$\mathcal{G}=\left( \mathcal{V},\mathcal{%
E}\right) $\textit{\ be a graph with }$n$\textit{\ vertices, }$n\geq 2.$%
\textit{\ Let }$k\in \left\langle n-1\right\rangle .$\textit{\ Then the
following statements are equivalent.}

(i)\textit{\ The graph }$\mathcal{G}$\textit{\ is }$k$\textit{-connected.}

(ii)\textit{\ The adjacency matrix of }$\mathcal{G}$\textit{\ is a }$%
g_{k}^{+}$\textit{-matrix.}

(iii)\textit{\ }$\forall \mathcal{Y},$\textit{\ }$\emptyset \neq \mathcal{Y}%
\subset \mathcal{V},$\textit{\ }$\exists \mathcal{X},$\textit{\ }$\emptyset
\neq \mathcal{X}\subseteq \mathcal{Y}^{c},$\textit{\ }$\left| \mathcal{X}%
\right| \geq \min \left( k,\left| \mathcal{Y}^{c}\right| \right) ,$\textit{\
and }$\forall V\in \mathcal{X},$\textit{\ }$\exists W\in \mathcal{Y}$\textit{%
\ such that }$\left[ V,W\right] \in \mathcal{E}.$

\smallskip

(iv) $\left| D_{\mathcal{Y}}\right| \geq \min \left( k,\left| \mathcal{Y}%
^{c}\right| \right) ,$ $\forall \mathcal{Y},$ $\emptyset \neq \mathcal{Y}%
\subset \mathcal{V}$.

\smallskip

\textit{Proof}. Definition 3.1, Theorem 3.2, and Theorem 3.25 and its proof. 
$\square $

\smallskip

We considered nondirected graphs. The directed graphs can also be
considered, the strongly connected directed graphs can be considered, ...
(For the directed graph (digraph) theory, see, \textit{e.g.}, [1].)

\bigskip

\begin{center}
\textbf{4. APPLICATIONS, I}
\end{center}

\bigskip

In this section, we give the first applications of $G^{+}$ method, more
exactly, of Theorems 2.6 and 2.8. The results, old and new results, are
mainly for the irreducible matrices and for the primitive ones, in some
cases the matrices being, moreover, $g_{k}^{+}$-matrices (not only
irreducible).

\smallskip

Recall that the row-allowable matrices and column-allowable ones were
defined in Section 2 (see Definition 2.1). The next theorem refers to these
matrices --- it is simple, but very useful.

\smallskip

\textbf{THEOREM 4.1.} \textit{Let }$P_{1}\in N_{n_{1},n_{2}}$\textit{\ and }$%
P_{2}\in N_{n_{2},n_{3}},$ $n_{1},$ $n_{2},$ $n_{3}\geq 1.$

(i) \textit{If }$P_{1}$\textit{\ is row-allowable and }$\left( P_{2}\right)
^{\left\{ j\right\} }>0$ (\textit{i.e., the column }$j$\textit{\ of }$P_{2}$%
\textit{\ is positive}), \textit{then\ }$\left( P_{1}P_{2}\right) ^{\left\{
j\right\} }>0$\textit{, where }$j\in \left\langle n_{3}\right\rangle .$

(ii) \textit{If }$\left( P_{1}\right) _{\left\{ i\right\} }>0$ (\textit{%
i.e., the row }$i$\textit{\ of }$P_{1}$\textit{\ is positive}) \textit{and }$%
P_{2}$\textit{\ is column-allowable}, \textit{then }$\left(
P_{1}P_{2}\right) _{\left\{ i\right\} }>0,$\textit{\ where} $i\in
\left\langle n_{1}\right\rangle $.

(iii) \textit{If }$\left( P_{1}\right) ^{\left\{ j\right\} }>0,$ \textit{%
where }$j\in \left\langle n_{2}\right\rangle ,$\textit{\ and }$P_{2}$\textit{%
\ is row-allowable, then\ }$\exists k\in \left\langle n_{3}\right\rangle $%
\textit{\ such that }$\left( P_{1}P_{2}\right) ^{\left\{ k\right\} }>0.$

(iv) \textit{If }$P_{1}$\textit{\ is column-allowable and }$\left(
P_{2}\right) _{\left\{ i\right\} }>0,$\textit{\ where }$i\in \left\langle
n_{2}\right\rangle $\textit{, then }$\exists k\in \left\langle
n_{1}\right\rangle $\textit{\ such that }$\left( P_{1}P_{2}\right) _{\left\{
k\right\} }>0$\textit{.}

(v) \textit{If }$P_{1}$\textit{\ is row-allowable and }$P_{2}$\textit{\ is
positive, then }$P_{1}P_{2}$\textit{\ is positive.}

(vi) \textit{If }$P_{1}$\textit{\ is positive and }$P_{2}$\textit{\ is
column-allowable, then }$P_{1}P_{2}$\textit{\ is positive.}

\smallskip

\textit{Proof}. (i) Obvious. We can even apply Theorem 2.8. Suppose that $%
P_{1}$ is row-allowable and $\left( P_{2}\right) ^{\left\{ j\right\} }>0.$
Then $P_{1}$ is sum-positive $\left\langle n_{1}\right\rangle \times
\left\langle n_{2}\right\rangle $ and $P_{2}$ is sum-positive $\left\langle
n_{2}\right\rangle \times \left\{ j\right\} .$ It follows from Theorem 2.8
that $\left( P_{1}P_{2}\right) ^{\left\{ j\right\} }>0$.

(ii) The transpose operation and (i).

(iii)-(vi) These are left to the reader --- for (v)-(vi), see also Remark
2.4. $\square $

\smallskip

Recall that the reducible matrices and irreducible ones were defined in
Section 3 (see Definitions 3.4 and 3.5). Recall that the $g_{k}^{+}$%
-matrices (in particular, the irreducible ones) are row-allowable (see
Theorem 3.10).

\smallskip

\textbf{THEOREM 4.2.} (See, \textit{e.g.}, [9, p. 507].) \textit{Let }$P\in
N_{n}$\textit{, }$n\geq 2.$\textit{\ Then }$P$\textit{\ is irreducible if
and only if} 
\[
\left( I+P\right) ^{n-1}>0 
\]
\noindent (\textit{recall that }$I$\textit{\ is the identity matrix})\textit{%
.}

\smallskip

\textit{Proof}. ``$\Longleftarrow $'' Obviously, $I+P$ is irreducible ($I+P$
is reducible $\Longrightarrow \left( I+P\right) ^{n-1}\ngtr 0$). So, $P$ is
irreducible.

``$\Longrightarrow $'' Set $Q=I+P.$ $Q$ is irreducible because $P$ is
irreducible and $Q\geq P.$ Moreover, we have $Q_{ii}>0,$ $\forall i\in
\left\langle n\right\rangle .$ We show that $Q^{n-1}>0,$ equivalently, that
(see Remark 2.4) 
\[
\left( Q^{n-1}\right) ^{\left\{ j\right\} }>0,\text{ }\forall j\in
\left\langle n\right\rangle . 
\]

Let $j\in \left\langle n\right\rangle .$

\smallskip

\textit{Case} 1. $Q^{\left\{ j\right\} }>0.$ Set $t=t\left( j\right) =1,$
and we have $\left( Q^{t}\right) ^{\left\{ j\right\} }>0.$

\smallskip

\textit{Case} 2. $Q^{\left\{ j\right\} }\ngtr 0.$ This case holds when $%
n\geq 3$ ($n=2,$ $Q_{ii}>0,$ $\forall i\in \left\langle n\right\rangle ,$
and $Q$ is irreducible $\Longrightarrow $ $Q>0$). Set --- a tail-to-head
construction --- 
\[
B_{1}=\left\{ j\right\} 
\]

\noindent and 
\[
B_{u+1}=B_{u}\cup \left\{ i\left| \text{ }i\in \left\langle n\right\rangle
-B_{u}\text{ and }\exists k\in B_{u}\text{ such that }Q_{ik}>0\right.
\right\} , 
\]

\noindent $\forall u=u\left( j\right) \geq 1$ with $B_{u}\subset
\left\langle n\right\rangle .$ By the definition of sets $B_{1}$ and $%
B_{u+1},$ $u\geq 1$ with $B_{u}\subset \left\langle n\right\rangle ,$ and
Theorem 3.7(iii) ($Q$ is irreducible...) we have 
\[
B_{u}\subset B_{u+1},\text{ }\forall u\geq 1\text{ with }B_{u}\subset
\left\langle n\right\rangle . 
\]
\noindent Since $Q^{\left\{ j\right\} }\ngtr 0,$ we have $B_{2}\subset
\left\langle n\right\rangle .$ Since $B_{u}\subset B_{u+1}\subseteq
\left\langle n\right\rangle ,$ $\forall u\geq 1$ with $B_{u}\subset
\left\langle n\right\rangle ,$ it follows that $\exists u_{0}=u_{0}\left(
j\right) \in \left\langle n\right\rangle -\left\{ 1,2\right\} $ such that $%
B_{u_{0}}=\left\langle n\right\rangle $ (recall that $n\geq 3$). Set $%
t=t\left( j\right) =u_{0}-1.$ Obviously, $t\in \left\langle n-1\right\rangle 
$. By the definition of sets $B_{1}$ and $B_{u+1},$ $u\geq 1$ with $%
B_{u}\subset \left\langle n\right\rangle ,$ and the fact that $Q_{ii}>0,$ $%
\forall i\in \left\langle n\right\rangle ,$ we have 
\[
B_{1}=\left\{ j\right\} \leftarrow B_{2}\leftarrow ...\leftarrow
B_{t}\leftarrow B_{t+1}=\left\langle n\right\rangle 
\]

\noindent ($u_{0}=t+1$). By Theorem 2.8, see also Remark 2.9, we have $%
\left( Q^{t}\right) ^{\left\{ j\right\} }>0$.

From Cases 1 and 2, we have $\left( Q^{t}\right) ^{\left\{ j\right\} }>0,$
where $t\in \left\langle n-1\right\rangle ,$%
\[
t=\left\{ 
\begin{array}{ll}
1\smallskip & \text{if }Q^{\left\{ j\right\} }>0, \\ 
u_{0}-1 & \text{if }Q^{\left\{ j\right\} }\ngtr 0.
\end{array}
\right. 
\]

\noindent If $t=n-1,$ no problem ($\left( Q^{n-1}\right) ^{\left\{ j\right\}
}>0$). If $1\leq t<n-1,$ by Theorem 4.1(i) we have 
\[
\left( Q^{n-1}\right) ^{\left\{ j\right\} }=\left( Q^{n-1-t}Q^{t}\right)
^{\left\{ j\right\} }=Q^{n-1-t}\left( Q^{t}\right) ^{\left\{ j\right\} }>0.%
\text{ }\square 
\]

\smallskip

Let $x\in \Bbb{R}$. Set $\left\lfloor x\right\rfloor =\max \left\{ k\left| 
\text{ }k\in \Bbb{Z}\text{ and }k\leq x\right. \right\} .$

\smallskip

The next result is a generalization of Theorem 4.2, ``$\Longrightarrow $''.

\smallskip

\textbf{THEOREM 4.3.} \textit{Let }$P\in N_{n},$\textit{\ }$n\geq 2$\textit{%
. Let }$k\in \left\langle n-1\right\rangle $\textit{. If }$P\in \overline{G}%
_{n,k}^{+}$\textit{, then} 
\[
\left( I+P\right) ^{m}>0, 
\]
\noindent \textit{where }$m=\left\lfloor \frac{n-2}{k}\right\rfloor +1.$

\smallskip

\textit{Proof}. Theorem 4.10 --- a more general result. $\square $

\smallskip

\textit{Remark} 4.4. (a) $\left( I+P\right) ^{m}>0$ is only a necessary
condition for a nonnegative $n\times n$ matrix $P$ be a $g_{k}^{+}$-matrix,
where $m=\left\lfloor \frac{n-2}{k}\right\rfloor +1.$ Indeed, considering 
\[
P=\left( 
\begin{array}{cccc}
\ast & * & * & * \\ 
\ast & * & 0 & 0 \\ 
0 & 0 & * & * \\ 
\ast & * & * & *
\end{array}
\right) , 
\]

\noindent where ``$*$'' stands for a positive entry, we have $\left(
I+P\right) ^{2}>0,$ $2=\left\lfloor \frac{4-2}{2}\right\rfloor +1$, but $P$
is not a $g_{2}^{+}$-matrix because (see Theorem 3.2) $D_{\left\{
3,4\right\} }=\left\{ 1\right\} ,$ $\left| D_{\left\{ 3,4\right\} }\right|
=1<\min \left( 2,2\right) =2.$

\smallskip

(b) By Theorem 4.3, $\left( I+P\right) ^{m}\ngtr 0\Longrightarrow P\notin 
\overline{G}_{n,k}^{+}$, where $P\in N_{n},$... So, we have a method to show
that a matrix $P\in N_{n}$ is not a $g_{k}^{+}$-matrix. In particular, due
to Theorem 3.25, we have a method to show that a graph with $n$ vertices, $%
n\geq 2,$ is not $k$-connected.

\smallskip

\textit{Definition} 4.5. (See, \textit{e.g.}, [9, p. 516] and [14, p. 47].)
Let $P\in N_{n},$ $n\geq 1$. We say that $P$ is \textit{primitive} if it is
irreducible and has only one eigenvalue of maximum modulus.

\smallskip

\textbf{THEOREM 4.6.} (See, \textit{e.g.}, [9, p. 516] and [14, p. 49].) 
\textit{Let }$P\in N_{n},$\textit{\ }$n\geq 1$\textit{. Then }$P$\textit{\
is primitive if and only if }$\exists m\geq 1$\textit{\ }($m\in \Bbb{N}$)%
\textit{\ such that }$P^{m}>0.$

\smallskip

\textit{Proof}. See, \textit{e.g.}, [9, pp. 516$-$517] and [14, pp. 49$-$%
50]. $\square $

\smallskip

Let $P\in N_{n}$ be a primitive matrix. Set (see, \textit{e.g.}, [9, p. 519]
and [21, p. 56]) 
\[
\gamma \left( P\right) =\text{the least }k\geq 1\text{ (}k\in \Bbb{N}\text{)
such that }P^{k}>0. 
\]

\noindent $\gamma \left( P\right) $ is called the \textit{index of
primitivity of }$P.$

\smallskip

\textit{Remark} 4.7. In general, it is not easy to compute $\gamma \left(
P\right) .$ In some cases (of interest or not) we can obtain good upper
bounds for $\gamma \left( P\right) $ --- the next result is for such a case.

\smallskip

\textbf{THEOREM 4.8.} (See, \textit{e.g.}, [9, p. 517]). \textit{Let }$P\in
N_{n},$\textit{\ }$n\geq 2$\textit{, be an irreducible matrix. If all the
main diagonal entries of }$P$\textit{\ are positive, then} 
\[
P^{n-1}>0 
\]
\noindent (\textit{and, therefore, }$P$\textit{\ is primitive and }$\gamma
\left( P\right) \leq n-1$)\textit{.}

\smallskip

\textit{Proof}. Similar to the proof of Theorem 4.2, ``$\Longrightarrow $''. 
$\square $

\smallskip

\textit{Remark} 4.9. Theorem 4.8 can be generalized as follows. Let $P\in
N_{n},$\ $n\geq 2$, be an irreducible matrix. Suppose that all the main
diagonal entries of $P$\ are positive. Let $P_{1},P_{2},...,P_{n-1}\in
N_{n}. $ Suppose that $P_{1},P_{2},...,P_{n-1}$ have the pattern $P$ (see
Definition 3.17). Then $P_{1}P_{2}...P_{n-1}>0.$ The proof is left to the
reader.

\smallskip

In this article, for simplification, the generalizations as that from Remark
4.9 will be omitted. Another generalization of Theorem 4.8 is given in the
next result.

\smallskip

\textbf{THEOREM 4.10.}\textit{\ Let }$P_{1},$ $P_{2},$ $...,$ $P_{m}\in 
\overline{G}_{n,k}^{+},$\textit{\ where }$n\geq 2$\textit{, }$k\in
\left\langle n-1\right\rangle $\textit{, and }$m=\left\lfloor \frac{n-2}{k}%
\right\rfloor +1.$\textit{\ Suppose that all the main diagonal entries of }$%
P_{1},$ $P_{2},$ $...,$ $P_{m}$\textit{\ are positive. Then} 
\[
P_{1}P_{2}...P_{m}>0. 
\]
\noindent \textit{In particular, if }$P_{1}=P_{2}=...=P_{m}:=P,$\textit{\
then }$P^{m}>0$\textit{\ and, as a result, }$P$\textit{\ is primitive and} 
\[
\gamma \left( P\right) \leq m=\left\lfloor \frac{n-2}{k}\right\rfloor +1. 
\]

\smallskip

\textit{Proof}. The proof is somehow similar to those of Theorems 4.2, ``$%
\Longrightarrow $'', and 4.8. We show that 
\[
\left( P_{1}P_{2}...P_{m}\right) ^{\left\{ j\right\} }>0,\text{ }\forall
j\in \left\langle n\right\rangle . 
\]

Let $j\in \left\langle n\right\rangle .$

\smallskip

\textit{Case} 1. $\left( P_{m}\right) ^{\left\{ j\right\} }>0.$ Set $%
t=t\left( j\right) =0.$

\smallskip

\textit{Case} 2. $\left( P_{m}\right) ^{\left\{ j\right\} }\ngtr 0.$ This
case holds when $n\geq 3$ and $m\geq 2$ ($n=2,$ $\left( P_{m}\right)
_{ii}>0, $ $\forall i\in \left\langle n\right\rangle ,$ and $P_{m}$ is
irreducible $\Longrightarrow $ $P_{m}>0$ ($P_{m}$ is irreducible because $%
P_{m}\in \overline{G}_{n,k}^{+}$); $P_{m}\in \overline{G}_{n,k}^{+},$ $%
\left( P_{m}\right) _{jj}>0,$ and $k=n-1$ $\Longrightarrow $ $\left(
P_{m}\right) ^{\left\{ j\right\} }>0$ (more generally, $P_{m}\in \overline{G}%
_{n,k}^{+},$ $\left( P_{m}\right) _{ii}>0,$ $\forall i\in \left\langle
n\right\rangle ,$ and $k=n-1$ $\Longrightarrow $ $P_{m}>0$); so, $k\in
\left\langle n-2\right\rangle ,$ and, further, we obtain $m\geq 2$). Set ---
a tail-to-head construction too --- 
\[
B_{1}=\left\{ j\right\} 
\]

\noindent and 
\[
B_{u+1}=B_{u}\cup C_{u},\text{ }\forall u=u\left( j\right) \geq 1\text{ with 
}B_{u}\subset \left\langle n\right\rangle , 
\]
\noindent where 
\[
C_{u}=\left\{ i\left| \text{ }i\in \left\langle n\right\rangle -B_{u}\text{
and }\exists k\in B_{u}\text{ such that }\left( P_{m-u+1}\right)
_{ik}>0\right. \right\} , 
\]

\noindent $\forall u\geq 1$ with $B_{u}\subset \left\langle n\right\rangle .$
Obviously, $B_{u}\cap C_{u}=\emptyset ,$ $\forall u\geq 1$ with $%
B_{u}\subset \left\langle n\right\rangle .$ Since $P_{l}\in \overline{G}%
_{n,k}^{+},$ $\forall l\in \left\langle m\right\rangle $, on the one hand,
using the definition of sets $B_{1}$ and $B_{u+1},$ $u\geq 1$ with $%
B_{u}\subset \left\langle n\right\rangle ,$ Remark 3.3(b), and Theorems
3.7(iii) and 3.9, we have 
\[
B_{u}\subset B_{u+1}=B_{u}\cup C_{u},\text{ }\forall u\geq 1\text{ with }%
B_{u}\subset \left\langle n\right\rangle , 
\]
\noindent and, on the other hand, we have 
\[
\left| B_{u+1}\right| =\left| B_{u}\right| +\left| C_{u}\right| \geq \left\{ 
\begin{array}{ll}
\left| B_{u}\right| +k\smallskip & \text{if }B_{u+1}\subset \left\langle
n\right\rangle , \\ 
\left| B_{u}\right| +\left| \left\langle n\right\rangle -B_{u}\right| & 
\text{if }B_{u+1}=\left\langle n\right\rangle ,
\end{array}
\right. 
\]
\[
\geq \left\{ 
\begin{array}{ll}
\left| B_{u-1}\right| +2k\smallskip & \text{if }B_{u+1}\subset \left\langle
n\right\rangle , \\ 
\left| B_{u-1}\right| +k+\left| \left\langle n\right\rangle -B_{u}\right| & 
\text{if }B_{u+1}=\left\langle n\right\rangle ,
\end{array}
\right. 
\]
\[
\geq ...\geq \left\{ 
\begin{array}{ll}
1+uk\smallskip & \text{if }B_{u+1}\subset \left\langle n\right\rangle , \\ 
1+\left( u-1\right) k+\left| \left\langle n\right\rangle -B_{u}\right| & 
\text{if }B_{u+1}=\left\langle n\right\rangle ,
\end{array}
\right. 
\]
$\forall u\geq 1$ with $B_{u}\subset \left\langle n\right\rangle .$ Set 
\[
a=\left| \left\langle n\right\rangle -B_{u}\right| \text{ if }%
B_{u+1}=\left\langle n\right\rangle ; 
\]
\noindent obviously, $1\leq a\leq k$ (see the definition of $g_{k}^{+}$%
-matrices again). If $B_{u+1}=\left\langle n\right\rangle $, we have $\left|
B_{u+1}\right| =n,$ so, in this case, 
\[
n\geq 1+\left( u-1\right) k+a. 
\]
\noindent Further, we obtain 
\[
u\leq \frac{n-1+k-a}{k}=\frac{n-1}{k}+\frac{k-a}{k}\leq 
\]
\[
\leq \frac{n-1}{k}+\frac{k-1}{k}=\frac{n-2}{k}+1. 
\]
\noindent Therefore, 
\[
u\leq \left\lfloor \frac{n-2}{k}\right\rfloor +1=m, 
\]
\noindent and, further, 
\[
u+1\leq m+1, 
\]
\noindent and thus we justified that the number of matrices we need must be $%
m$ (not $m+1$). Since $\left( P_{m}\right) ^{\left\{ j\right\} }\ngtr 0,$ we
have $B_{2}\subset \left\langle n\right\rangle .$ Since $B_{u}\subset
B_{u+1}\subseteq \left\langle n\right\rangle ,$ $\forall u\geq 1$ with $%
B_{u}\subset \left\langle n\right\rangle ,$ it follows that $\exists
u_{0}=u_{0}\left( j\right) \in \left\langle m+1\right\rangle -\left\{
1,2\right\} $ such that $B_{u_{0}}=\left\langle n\right\rangle $ (see above (%
$n\geq 3,$ $m\geq 2,$ ...)). Set $t=t\left( j\right) =u_{0}-2.$ Obviously, $%
t\in \left\langle m-1\right\rangle $. By the definition of sets $B_{1}$ and $%
B_{u+1},$ $u\geq 1$ with $B_{u}\subset \left\langle n\right\rangle ,$ and
the fact that $\left( P_{l}\right) _{ii}>0,$ $\forall l\in \left\langle
m\right\rangle ,\forall i\in \left\langle n\right\rangle $ we have 
\[
B_{1}=\left\{ j\right\} \leftarrow B_{2}=B_{1}\cup C_{1}\leftarrow
...\leftarrow B_{t+1}=B_{t}\cup C_{t}\leftarrow B_{t+2}=\left\langle
n\right\rangle 
\]

\noindent ($u_{0}=t+2$). By Theorem 2.8, see also Remark 2.9, we have 
\[
\left( P_{m-t}P_{m-t+1}...P_{m}\right) ^{\left\{ j\right\} }>0 
\]

\noindent (recall that $t\in \left\langle m-1\right\rangle $ --- so, $t\geq
1 $).

From Cases 1 and 2, we have $\left( P_{m-t}P_{m-t+1}...P_{m}\right)
^{\left\{ j\right\} }>0,$ where $t\in \left\langle \left\langle
m-1\right\rangle \right\rangle ,$%
\[
t=\left\{ 
\begin{array}{ll}
0\smallskip & \text{if }\left( P_{m}\right) ^{\left\{ j\right\} }>0, \\ 
u_{0}-2 & \text{if }\left( P_{m}\right) ^{\left\{ j\right\} }\ngtr 0.
\end{array}
\right. 
\]

\noindent If $t=m-1,$ no problem ($\left( P_{1}P_{2}...P_{m}\right)
^{\left\{ j\right\} }>0$). If $0\leq t<m-1,$ by Theorem 4.1(i) we have 
\[
\left( P_{1}P_{2}...P_{m}\right) ^{\left\{ j\right\} }=\left( \left(
P_{1}P_{2}...P_{m-t-1}\right) \left( P_{m-t}P_{m-t+1}...P_{m}\right) \right)
^{\left\{ j\right\} }= 
\]
\[
=\left( P_{1}P_{2}...P_{m-t-1}\right) \left( P_{m-t}P_{m-t+1}...P_{m}\right)
^{\left\{ j\right\} }>0.\text{ }\square 
\]

\smallskip

\textit{Example} 4.11. Let $P\in N_{9},$%
\[
P=\left( 
\begin{array}{ccccccccc}
\ast & 0 & 0 & * & * & * & 0 & 0 & 0 \\ 
0 & * & 0 & * & * & * & 0 & 0 & 0 \\ 
0 & 0 & * & * & * & * & 0 & 0 & 0 \\ 
0 & 0 & 0 & * & 0 & 0 & * & * & * \\ 
0 & 0 & 0 & 0 & * & 0 & * & * & * \\ 
0 & 0 & 0 & 0 & 0 & * & * & * & * \\ 
\ast & * & * & 0 & 0 & 0 & * & 0 & 0 \\ 
\ast & * & * & 0 & 0 & 0 & 0 & * & 0 \\ 
\ast & * & * & 0 & 0 & 0 & 0 & 0 & *
\end{array}
\right) , 
\]

\noindent where ``$*$'' stands for a positive entry. $P\in \overline{G}%
_{9,3}^{+}.$ (See Examples 3.15 and 3.18 again.) By Theorem 4.8 we have $%
P^{8}>0$ while by Theorem 4.10 we have $P^{3}>0.$ By direct computation, $%
P^{2}\ngtr 0$ and $P^{3}>0.$

\smallskip

\textit{Remark} 4.12. Theorem 4.10 can also be generalized. Let $P_{1}\in 
\overline{G}_{n,k_{1}}^{+},$ $P_{2}\in \overline{G}_{n,k_{2}}^{+},$ $...,$ $%
P_{m}\in \overline{G}_{n,k_{m}}^{+},$ where... --- the study of this case is
left to the reader.

\smallskip

We will give another generalization, a known one, of Theorem 4.8. To prove
this generalization using our approach, we need the following result.

\smallskip

\textbf{THEOREM 4.13.}\textit{\ Let }$P\in N_{n},$\textit{\ }$n\geq 2$%
\textit{, be an irreducible matrix. Suppose that }$P$\textit{\ has }$d$%
\textit{\ positive main diagonal entries, }$d\in \left\langle n\right\rangle 
$\textit{\ --- suppose that }$P_{i_{u}i_{u}}>0,$ $\forall u\in \left\langle
d\right\rangle ,$\textit{\ where }$i_{1},i_{2},...,i_{d}\in \left\langle
n\right\rangle ,$ $i_{u}\neq i_{v},$ $\forall u,v,$ $u\neq v.$\textit{\ Let }%
$W=\left\{ i_{1},i_{2},...,i_{d}\right\} .$\textit{\ Then} 
\[
\left( P^{n-d}\right) ^{W} 
\]
\noindent \textit{is row-allowable }(\textit{sum-positive on }$\left\langle
n\right\rangle \times W$) \textit{and} 
\[
\left( P^{n-d}\right) _{W} 
\]
\noindent \textit{is column-allowable.}

\smallskip

\textit{Proof}. \textit{Case} 1. $d=n.$ Obvious.

\textit{Case} 2. $1\leq d<n.$ First, we show that $\left( P^{n-d}\right)
^{W} $ is row-allowable. Since $P_{i_{u}i_{u}}>0,$ $\forall u\in
\left\langle d\right\rangle ,$ it follows that $\left( P^{n-d}\right)
_{W}^{W}$ is row-allowable. It remains to show that $\left( P^{n-d}\right)
_{W^{c}}^{W}$ is row-allowable ($W^{c}$ is the complement $W$).

\smallskip

Let $i\in W^{c}.$ We show that $P_{ij}^{n-d}>0$ for some $j\in W.$ We can
either have $P_{ij}>0$ for some $j\in W$ or $P_{ij}=0,$ $\forall j\in W$. If 
$P_{ij}>0,$ we have $P\in \overline{G}_{\left\{ i\right\} ,\left\{ j\right\}
}^{+}$. Obviously, $P^{t}\in \overline{G}_{\left\{ k\right\} ,\left\{
k\right\} }^{+},$ $\forall t\geq 0,$ $\forall k\in W.$ By Theorem 2.6(i), $%
P^{n-d}\in \overline{G}_{\left\{ i\right\} ,\left\{ j\right\} }^{+}$ ($%
P^{n-d}=PP^{n-d-1},$ ...), so, $P_{ij}^{n-d}>0.$ If $P_{ij}=0,$ $\forall
j\in W$ --- this subcase holds only when $\left| W^{c}\right| >1$ because $P$
is irreducible ---, then $\exists l\in \left\langle n-d-1\right\rangle ,$ $%
\exists j_{1},j_{2},...,j_{l}\in W^{c},$ $j_{u}\neq j_{v},$ $\forall u,v\in
\left\langle l\right\rangle ,$ $u\neq v,$ and $j_{u}\neq i,$ $\forall u\in
\left\langle l\right\rangle ,$ $\exists j_{l+1}\in W$ such that $P_{ij_{1}},$
$P_{j_{1}j_{2}},$ $...,$ $P_{j_{l}j_{l+1}}>0$ (because $P$ is irreducible
and $\left| W^{c}-\left\{ i\right\} \right| =n-d-1>0$). Further, we have $%
P\in \overline{G}_{\left\{ i\right\} ,\left\{ j_{1}\right\} }^{+}\cap 
\overline{G}_{\left\{ j_{1}\right\} ,\left\{ j_{2}\right\} }^{+}\cap ...\cap 
\overline{G}_{\left\{ j_{l}\right\} ,\left\{ j_{l+1}\right\} }^{+},$ so, by
Theorem 2.6(ii), $P^{l+1}\in \overline{G}_{\left\{ i\right\} ,\left\{
j_{l+1}\right\} }^{+}.$ Therefore, $P_{ij_{l+1}}^{l+1}>0.$

If $l=n-d-1,$ no problem --- we have $P_{ij}^{n-d}>0$ for some $j\in W$ ($%
j=j_{n-d}$).

If $1\leq l<n-d-1,$ using Theorem 2.6(i) for $P^{l+1}\in \overline{G}%
_{\left\{ i\right\} ,\left\{ j_{l+1}\right\} }^{+}$ and $P^{n-d-l-1}\in 
\overline{G}_{\left\{ j_{l+1}\right\} ,\left\{ j_{l+1}\right\} }^{+}$
(recall that $P^{t}\in \overline{G}_{\left\{ k\right\} ,\left\{ k\right\}
}^{+},$ $\forall t\geq 1$ (even $\forall t\geq 0$)$,$ $\forall k\in W$), we
have $P^{n-d}\in \overline{G}_{\left\{ i\right\} ,\left\{ j_{l+1}\right\}
}^{+},$ so, $P_{ij_{l+1}}^{n-d}>0.$

We conclude that $\left( P^{n-d}\right) ^{W}$ is row-allowable.

The above result leads to the fact that $\left( \left( ^{\text{t}}P\right)
^{n-d}\right) ^{W}$ is row-allowable ($^{\text{t}}P$ is the transpose of $P$%
). So, $\left( P^{n-d}\right) _{W}$ is column-allowable. $\square $

\smallskip

Now, we give the generalization of Theorem 4.8 we promised.

\smallskip

\textbf{THEOREM 4.14.} (See, \textit{e.g.}, [9, pp. 520$-$521].) \textit{%
Under the same conditions as in Theorem }4.13\textit{\ we have} 
\[
P^{2n-d-1}>0 
\]
\noindent (\textit{so}, $\gamma \left( P\right) \leq 2n-d-1$).

\smallskip

\textit{Proof}. \textit{Case} 1. $d=n.$ Theorem 4.8.

\textit{Case} 2. $d\in \left\langle n-1\right\rangle .$ Proceeding as in the
proof of Theorem 4.2, ``$\Longrightarrow $'' (for the columns $%
i_{1},i_{2},...,i_{d}$), we have 
\[
\left( P^{n-1}\right) ^{W}>0\text{ and }\left( \left( ^{\text{t}}P\right)
^{n-1}\right) ^{W}>0. 
\]
\noindent The latter result leads to $\left( P^{n-1}\right) _{W}>0.$ By
Theorem 4.13, $\left( P^{n-d}\right) ^{W}$ is row-allowable and $\left(
P^{n-d}\right) _{W}$ is column-allowable. Finally, by Theorem 4.1(v) we have 
\[
P^{2n-d-1}=P^{n-d}P^{n-1}\geq \left( P^{n-d}\right) ^{W}\left(
P^{n-1}\right) _{W}>0. 
\]

\noindent (The inequality $P^{2n-d-1}>0$ can also be obtained using Theorem
2.8 ($P^{n-d}\in \overline{G}_{\left\langle n\right\rangle ,W}^{+},$ $%
P^{n-1}\in \overline{G}_{W,\left\{ j\right\} }^{+},$ $\forall j\in
\left\langle n\right\rangle ,$ ...) or the fact that $\left( P^{n-1}\right)
^{W}>0$ and $\left( P^{n-d}\right) _{W}$ is column-allowable --- in the
latter case, by Theorem 4.1(vi) we have 
\[
P^{2n-d-1}=P^{n-1}P^{n-d}\geq \left( P^{n-1}\right) ^{W}\left(
P^{n-d}\right) _{W}>0.\text{) }\square 
\]

\smallskip

Theorem 4.14 can be generalized --- the next result is a generalization both
of Theorem 4.14 and of Theorem 4.10.

\smallskip

\textbf{THEOREM 4.15.}\textit{\ Let }$P_{1},P_{2},...,P_{m+1}\in N_{n}$%
\textit{\ be irreducible matrices, where }$m=\left\lfloor \frac{n-2}{k}%
\right\rfloor +1,$ $n\geq 2$\textit{, and }$k\in \left\langle
n-1\right\rangle $\textit{. Let }$W=\left\{ i_{1},i_{2},...,i_{d}\right\} ,$ 
$d\in \left\langle n\right\rangle .$

(i) \textit{If }$\left( P_{1}\right) ^{W}$\textit{\ is row-allowable }(%
\textit{irreducible or not})\textit{, }$\left( P_{l}\right) _{i_{u}i_{u}}>0,$
$\forall l\in \left\langle m+1\right\rangle -\left\{ 1\right\} ,$ $\forall
u\in \left\langle d\right\rangle ,$ \textit{and }$P_{2},P_{3},...,P_{m+1}\in 
\overline{G}_{n,k}^{+},$\textit{\ then} 
\[
P_{1}P_{2}...P_{m+1}>0. 
\]
\noindent \textit{This result holds, in particular, for }$%
P_{2}=P_{3}=...=P_{m+1}:=P$ and $P_{1}=P^{n-d},$\textit{\ and we have} 
\[
P^{n+m-d}>0 
\]
\noindent (\textit{so, }$\gamma \left( P\right) \leq n+\left\lfloor \frac{n-2%
}{k}\right\rfloor +1-d$).

(ii) \textit{If }$\left( P_{l}\right) _{i_{u}i_{u}}>0,$ $\forall l\in
\left\langle m\right\rangle ,$ $\forall u\in \left\langle d\right\rangle ,$%
\textit{\ }$P_{1},P_{2},...,P_{m}\in \overline{G}_{n,k}^{+},$\textit{\ and }$%
\left( P_{m+1}\right) _{W}$\textit{\ is column-allowable }(\textit{%
irreducible or not})\textit{, then} 
\[
P_{1}P_{2}...P_{m+1}>0. 
\]
\noindent \textit{This result holds, in particular, for }$%
P_{1}=P_{2}=...=P_{m}:=P$ and $P_{m+1}=P^{n-d},$\textit{\ and we have} 
\[
P^{n+m-d}>0. 
\]

\smallskip

\textit{Proof}. (i) \textit{Case} 1. $d=n.$ Theorems 4.1(v) and 4.10.

\textit{Case} 2. $d\in \left\langle n-1\right\rangle .$ Proceeding as in the
proof of Theorem 4.10 for the matrices $^{\text{t}}P_{m+1},$ $^{\text{t}%
}P_{m},$ $...,$ $^{\text{t}}P_{2}$ and columns $i_{1},i_{2},...,i_{d}$, we
have 
\[
\left( ^{\text{t}}P_{m+1}{}^{\text{t}}P_{m}...^{\text{t}}P_{2}\right)
^{W}>0. 
\]
\noindent Further, we have $\left( ^{\text{t}}\left(
P_{2}P_{3}...P_{m+1}\right) \right) ^{W}>0.$ So, $\left(
P_{2}P_{3}...P_{m+1}\right) _{W}>0.$ Finally, by Theorem 4.1(v) we have 
\[
P_{1}P_{2}...P_{m+1}=P_{1}\left( P_{2}P_{3}...P_{m+1}\right) \geq \left(
P_{1}\right) ^{W}\left( P_{2}P_{3}...P_{m+1}\right) _{W}>0. 
\]

\smallskip

(ii) Similar to (i) (using Theorem 4.1(vi) instead of Theorem 4.1(v)). $%
\square $

\smallskip

\textit{Remark} 4.16. Theorem 4.10 is a special case of Theorem 4.15 --- it
is the case when $d=n$ (equivalently, $W=\left\langle n\right\rangle $), $%
P_{1}=I$ at (i), and $P_{m+1}=I$ at (ii).

\smallskip

\textit{Definition} 4.17. (See, \textit{e.g.}, [21, p. 11]; see, \textit{e.g.%
}, also [9, p. 357].) Let $P\in N_{n},$\ $n\geq 1$. Let $i,$ $j\in
\left\langle n\right\rangle .$ Let $i_{0},i_{1},...,i_{u}\in \left\langle
n\right\rangle ,$ $u\geq 1,$ $i_{0}=i,$ $i_{u}=j.$ We say that $\left(
i_{0},i_{1},...,i_{u}\right) $ is a \textit{path/chain} (\textit{of} $P$) 
\textit{from }$i$\textit{\ to }$j$ if $P_{i_{0}i_{1}},$ $P_{i_{1}i_{2}},$ $%
...,$ $P_{i_{u-1}i_{u}}>0.$ $u$ is called the \textit{length of path from }$%
i $\textit{\ to }$j.$ If $i=j,$ the path is called a \textit{cycle from }$i$%
\textit{\ to }$i$\textit{\ }(or \textit{from }$i$\textit{\ to itself}).

\smallskip

\textit{Definition} 4.18. (See, \textit{e.g.}, [21, p. 16].) Let $P\in
N_{n}, $\ $n\geq 1$. Let $i\in \left\langle n\right\rangle .$ Suppose that
there exists a \textit{cycle from }$i$\textit{\ to }$i$. Set $d(i)=$ the
greatest common divisor of those $k\geq 1$ for which $P_{ii}^{k}>0$ ($%
P_{ii}^{k}=\left( P^{k}\right) _{ii}$). $d(i)$ is called \textit{period of }$%
i.$

\smallskip

\textbf{THEOREM 4.19.} (See, \textit{e.g.}, [21, pp. 17$-$18].) \textit{Let }%
$P\in N_{n},$\textit{\ }$n\geq 1,$\textit{\ be an irreducible matrix. Then }$%
\forall i,$\textit{\ }$j\in \left\langle n\right\rangle $\textit{\ we have }$%
d(i)=d(j).$

\smallskip

\textit{Proof}. See, \textit{e.g.}, [21, p. 17, Lemma 1.2, and p. 18,
Definition 1.6]. $\square $

\smallskip

\textit{Definition} 4.20. (See, \textit{e.g.}, [21, p. 18].) Let $P\in
N_{n}, $\ $n\geq 1,$ be an irreducible matrix. Set $d=d(1)$ (by Theorem
4.19, $d=d(1)=d(2)=...=d(n)$). We say that $P$ is \textit{periodic/cyclic }(%
\textit{with period }$d$) if $d>1.$ If $d=1$, we say that $P$ is \textit{%
aperiodic/acyclic}.

\smallskip

\textbf{THEOREM 4.21.} (See, \textit{e.g.}, [21, p. 21].) \textit{Let }$P\in
N_{n},$\textit{\ }$n\geq 1$\textit{. Then }$P$\textit{\ is primitive }(%
\textit{see Definition }4.5;\textit{\ see also Theorem }4.6)\textit{\ if and
only if it is aperiodic irreducible.}

\smallskip

\textit{Proof}. See, \textit{e.g.}, [21, p. 21]. $\square $

\smallskip

\textbf{THEOREM 4.22.} (See, \textit{e.g.}, [9, pp. 519$-$520].) \textit{Let 
}$P\in N_{n},$\textit{\ }$n\geq 1$\textit{\ be a primitive matrix. Let }$%
s\geq 1$\textit{\ be the smallest }(\textit{natural}) \textit{number for
which there exists a cycle }(\textit{of }$P$) \textit{with length }$s.$%
\textit{\ Then} 
\[
P>0\hspace{0.2cm}\text{\textit{if} }n=1 
\]
\noindent \textit{and} 
\[
P^{n+s\left( n-2\right) }>0\text{ (\textit{so}, }\gamma \left( P\right) \leq
n+s\left( n-2\right) \text{)}\hspace{0.2cm}\text{\textit{if} }n\geq 2. 
\]

\smallskip

\textit{Proof}. No problem when $n=1.$ Further, we consider that $n\geq 2,$
and show that 
\[
\left( P^{n+s\left( n-2\right) }\right) ^{\left\{ j\right\} }>0,\text{ }%
\forall j\in \left\langle n\right\rangle . 
\]

Fix $j\in \left\langle n\right\rangle .$ Fix a cycle of $P$ with length $s.$
Suppose that the cycle is $\left( i_{s},i_{s-1},...,i_{0}\right) $ with $%
i_{s}=i_{0}$, and, as a result, we have $P_{i_{s}i_{s-1}},$ $%
P_{i_{s-1}i_{s-2}},$ $...,$ $P_{i_{1}i_{0}}>0.$

\smallskip

\textit{Case} 1. $j$ belongs to the (fixed) cycle. Suppose that (no problem) 
$i_{s}=i_{0}=j.$ Set 
\[
U_{0}=\left\{ j\right\} \text{ and }U_{t+1}=\left\{ i\left| \text{ }i\in
\left\langle n\right\rangle \text{ and }\exists k\in U_{t}\text{ such that }%
P_{ik}>0\right. \right\} ,\forall t\geq 0. 
\]
\noindent Obviously, $i_{0}\in U_{ks},$ $i_{1}\in U_{ks+1},$ $...,$ $%
i_{s-1}\in U_{ks+s-1},$ $\forall k\geq 0.$ Obviously, we have 
\[
U_{0}=\left\{ j\right\} \leftarrow U_{1}\leftarrow ...\leftarrow
U_{s-1}\leftarrow 
\]
\[
\leftarrow U_{s}=U_{0}\cup U_{0}^{\left( 1\right) }\leftarrow
U_{s+1}=U_{1}\cup U_{1}^{\left( 1\right) }\leftarrow ...\leftarrow
U_{s+s-1}=U_{s-1}\cup U_{s-1}^{\left( 1\right) }\leftarrow 
\]
\[
\leftarrow U_{2s}=\!U_{s}\cup U_{0}^{\left( 2\right) }\leftarrow
U_{2s+1}=\!U_{s+1}\cup U_{1}^{\left( 2\right) }\leftarrow ...\leftarrow
U_{2s+s-1}=\!U_{s+s-1}\cup U_{s-1}^{\left( 2\right) }\leftarrow ..., 
\]
\noindent where $U_{0}^{\left( 1\right) }=U_{s}-U_{0},$ $U_{1}^{\left(
1\right) }=U_{s+1}-U_{1},...$ It follows that 
\[
U_{0}\subseteq U_{s}\subseteq U_{2s}\subseteq ..., 
\]
\[
U_{1}\subseteq U_{s+1}\subseteq U_{2s+1}\subseteq ..., 
\]
\[
\vdots 
\]
\[
U_{s-1}\subseteq U_{2s-1}\subseteq U_{3s-1}\subseteq ... 
\]
\noindent We cannot have 
\[
U_{ts}=U_{\left( t+1\right) s}\text{ if }U_{ts}\subset \left\langle
n\right\rangle , 
\]
\noindent where $t\geq 0,$ because $P,$ being primitive, is not
cyclic/periodic (see Theorem 4.21). More generally, we cannot have 
\[
U_{ts+w}=U_{\left( t+1\right) s+w}\text{ if }U_{ts+w}\subset \left\langle
n\right\rangle , 
\]
\noindent where $t\geq 0$ and $w\in \left\langle \left\langle
s-1\right\rangle \right\rangle .$

$U_{0}$ ($U_{0}=U_{0\cdot s}$) has one element. Since $P$ is aperiodic
irreducible (not periodic irreducible), it follows that $U_{s}$ ($%
U_{s}=U_{1\cdot s}$) has at least $2$ elements, $U_{2s}$ has at least $3$
elements, ..., $U_{\left( n-1\right) s}$ has at least $n$ elements ($n-1+1=n$%
), and, therefore, $U_{\left( n-1\right) s}=\left\langle n\right\rangle .$
Using Theorem 2.8 for $U_{0}=\left\{ j\right\} \leftarrow U_{1}\leftarrow
...\leftarrow U_{\left( n-1\right) s}=\left\langle n\right\rangle ,$ we have 
\[
\left( P^{\left( n-1\right) s}\right) ^{\left\{ j\right\} }>0. 
\]

\textit{Case} 2. $j$ does not belong to the cycle. In this case, there
exists a path from $i_{v}$ for some $v\in \left\langle \left\langle
s-1\right\rangle \right\rangle $ to $j$ with length at most $n-s$ ($s<n$
because $P$ is primitive, so, $n-s\geq 1$). Consider that this path is $%
\left( j_{0},j_{1},...,j_{z}\right) ,$ where $j_{0}=i_{v}$ and $j_{z}=j,$ $%
z\in \left\langle n-s\right\rangle .$ We have 
\[
\left\{ j\right\} =\left\{ j_{z}\right\} \leftarrow \left\{ j_{z-1}\right\}
\leftarrow ...\leftarrow \left\{ j_{1}\right\} \leftarrow \left\{
j_{0}\right\} =\left\{ i_{v}\right\} . 
\]
\noindent Further, for $i_{v},$ we use Case 1, and keeping notation for
sets, $U_{0},U_{1},$ $...,$ but, here, $U_{0}=\left\{ i_{v}\right\} ,$ and,
obviously, keeping the definitions for $U_{t+1},$ $t\geq 0,$ we have $%
U_{0}=\left\{ i_{v}\right\} \leftarrow U_{1}\leftarrow ...\leftarrow
U_{\left( n-1\right) s}=\left\langle n\right\rangle .$ Using Theorem 2.8 for 
\[
\left\{ j\right\} =\!\left\{ j_{z}\right\} \leftarrow \left\{
j_{z-1}\right\} \leftarrow ...\leftarrow \left\{ j_{1}\right\} \leftarrow
U_{0}=\!\left\{ j_{0}\right\} =\!\left\{ i_{v}\right\} \leftarrow
...\leftarrow U_{\left( n-1\right) s}=\!\left\langle n\right\rangle \!, 
\]
\noindent we have 
\[
\left( P^{z+\left( n-1\right) s}\right) ^{\left\{ j\right\} }>0. 
\]
\noindent Both when $z=n-s$ and when $1\leq z<n-s,$ we have 
\[
\left( P^{n-s+\left( n-1\right) s}\right) ^{\left\{ j\right\} }>0, 
\]
\noindent in the latter case using Theorem 4.1(i). Since $n-s+\left(
n-1\right) s=n+s\left( n-2\right) ,$ we have 
\[
\left( P^{n+s\left( n-2\right) }\right) ^{\left\{ j\right\} }>0. 
\]

From Cases 1 and 2, since $\max \left( \left( n-1\right) s,n+s\left(
n-2\right) \right) =n+s\left( n-2\right) ,$ by Theorem 4.1(i) we have 
\[
\left( P^{n+s\left( n-2\right) }\right) ^{\left\{ j\right\} }>0.\text{ }%
\square 
\]

\smallskip

Theorem 4.22 can be generalized.

\smallskip

\textbf{THEOREM 4.23.} \textit{Let }$P\in \overline{G}_{1,1}^{+}\cup 
\overline{G}_{n,k}^{+}$\textit{\ be a primitive matrix,\ where }$n\geq 2$%
\textit{\ and }$k\in \left\langle n-1\right\rangle .$ \textit{Let }$s\geq 1$%
\textit{\ be the smallest number for which there exists a cycle with length }%
$s.$\textit{\ Then} 
\[
P^{g}>0, 
\]
\noindent \textit{where} 
\[
g=\left\{ 
\begin{array}{ll}
1\medskip & \hspace{0.2cm}\text{\textit{if} }P\in \overline{G}_{1,1}^{+}, \\ 
\left\lfloor \frac{n-s-2}{k}\right\rfloor +2+s\left[ n-\max \left(
2,k\right) +1\right] & \hspace{0.2cm}\text{\textit{if} }P\in \overline{G}%
_{n,k}^{+}.
\end{array}
\right. 
\]

\noindent (\textit{For }$n\geq 2$ \textit{and} $k=1,$\textit{\ we have }$%
g=n+s\left( n-2\right) $\textit{\ --- as in Theorem} 4.22.)

\smallskip

\textit{Proof}. When $P\in \overline{G}_{1,1}^{+},$ we have $P\in N_{1},$ $%
P>0.$ So, no problem. Further, we consider that $P\in \overline{G}%
_{n,k}^{+}, $ $n\geq 2,$ $k\in \left\langle n-1\right\rangle ,$ and show
that 
\[
\left( P^{g}\right) ^{\left\{ j\right\} }>0,\text{ }\forall j\in
\left\langle n\right\rangle . 
\]

Fix $j\in \left\langle n\right\rangle .$ Fix a cycle of $P$ with length $s.$
Suppose that the cycle is $\left( i_{s},i_{s-1},...,i_{0}\right) $ with $%
i_{s}=i_{0}$ --- it follows that $P_{i_{s}i_{s-1}},$ $P_{i_{s-1}i_{s-2}},$ $%
...,$ $P_{i_{1}i_{0}}>0.$

\smallskip

\textit{Case} 1. $j$ belongs to the cycle. We proceed as in Case 1 from the
proof of Theorem 4.22, and have a little difference due to the fact that $%
P\in \overline{G}_{n,k}^{+}.$\textit{\ }Suppose that $i_{s}=i_{0}=j$ as
well, $U_{0}=\left\{ j\right\} $ as well, $U_{t+1}=...$ --- see there ---, $%
\forall t\geq 0,$ as well.

$U_{0}$ ($U_{0}=U_{0\cdot s}$) has one element. Since $P\in \overline{G}%
_{n,k}^{+}$ and is aperiodic irreducible, it follows that $U_{s}$ ($%
U_{s}=U_{1\cdot s}$) has at least $b$ elements, where $b=\max \left(
2,k\right) $ (not $b=k+1$; $b=2$ for $k=1$ and $k=2$; $b=k$ for $k\geq 3$), $%
U_{2s}$ has at least $b+1$ elements, $U_{3s}$ has at least $b+2$ elements,
..., $U_{\left( n-b+1\right) s}$ has at least $n$ elements ($b+n-b+1-1=n$),
and, therefore, $U_{\left( n-b+1\right) s}=\left\langle n\right\rangle .$
Using Theorem 2.8 for $U_{0}=\left\{ j\right\} \leftarrow U_{1}\leftarrow
...\leftarrow U_{\left( n-b+1\right) s}=\left\langle n\right\rangle ,$ we
have 
\[
\left( P^{\left( n-b+1\right) s}\right) ^{\left\{ j\right\} }>0. 
\]

\textit{Case} 2. $j$ does not belong to the cycle. Set 
\[
B_{1}=\left\{ j\right\} 
\]
\noindent and 
\[
B_{u}=\left\{ i\left| \text{ }i\in \left\langle n\right\rangle
-\bigcup_{t=1}^{u-1}B_{t}\text{ and }\exists k\in \bigcup_{t=1}^{u-1}B_{t}%
\text{ such that }P_{ik}>0\right. \right\} , 
\]

\noindent $\forall u\geq 2$ with $\left|
\bigcup\limits_{t=1}^{u-1}B_{t}\right| \leq n-s.$ Let $w$ be the smallest
(natural) number such that $\left| \bigcup\limits_{t=1}^{w}B_{t}\right|
>n-s. $ Obviously, $w\geq 2.$ Since $P\in \overline{G}_{n,k}^{+},$ we have 
\[
\left| B_{2}\right| \geq 1 
\]
\noindent when $w=2$ and 
\[
\left| B_{2}\right| \geq k,\text{ }\left| B_{3}\right| \geq k,\text{ }...,%
\text{ }\left| B_{w-1}\right| \geq k,\text{ }\left| B_{w}\right| \geq 1\text{
(}\left| B_{w}\right| \geq k\text{ or }1\leq \left| B_{w}\right| <k\text{)} 
\]
\noindent when $w\geq 3$. Further, since the sets $B_{1},B_{2},...,B_{w}$
are disjoint, we have 
\[
\left| \bigcup\limits_{t=1}^{w}B_{t}\right| =\left| B_{1}\right| +\left|
B_{2}\right| +...+\left| B_{w}\right| \geq 2+\left( w-2\right) k. 
\]
\noindent We must consider that 
\[
2+\left( w-2\right) k>n-s. 
\]
\noindent Further, we have 
\[
w>\frac{n-s-2}{k}+2, 
\]
\noindent so, 
\[
w=\left\lfloor \frac{n-s-2}{k}\right\rfloor +3. 
\]
\noindent By the definition of $w,$ $\exists l\in B_{w}$ such that $l$
belongs to the cycle. It follows that $l=i_{v}$ for some $v\in \left\langle
\left\langle s-1\right\rangle \right\rangle .$ We have either $B_{w}=\left\{
i_{v}\right\} $ or $B_{w}\supset \left\{ i_{v}\right\} .$ If $B_{w}\supset
\left\{ i_{v}\right\} ,$ we can work with $\left\{ i_{v}\right\} $ instead
of $B_{w}$ because if a matrix is sum-positive on $\left( C_{1}\cup
C_{2}\right) \times D,$ then it is sum-positive on $C_{1}\times D$ and on $%
C_{2}\times D.$ So, we can work with $\left\{ i_{v}\right\} $ in both cases.
We work with $\left\{ i_{v}\right\} $ in both cases, and, further, for $%
i_{v},$ we use Case 1, and keeping notation for sets, $U_{0},U_{1},$ $...,$
but, here, $U_{0}=\left\{ i_{v}\right\} ,$ and, obviously, keeping the
definition for $U_{t+1},$ $t\geq 0,$ we have 
\[
U_{0}=\left\{ i_{v}\right\} \leftarrow U_{1}\leftarrow ...\leftarrow
U_{\left( n-b+1\right) s}=\left\langle n\right\rangle . 
\]

\noindent By the definition of sets $B_{1},B_{2},...,B_{w}$ we have 
\[
\bigcup\limits_{t=1}^{u-1}B_{t}\leftarrow B_{u},\text{ }\forall u,\text{ }%
2\leq u\leq w. 
\]
\noindent But also we have 
\[
B_{u-1}\leftarrow B_{u},\text{ }\forall u,\text{ }2\leq u\leq w 
\]
\noindent --- we prove this statement. If $w=2$ (recall that $w\geq 2$), we
have only $B_{1}\leftarrow B_{2}$, so, no problem. Further, we consider that 
$w\geq 3.$ If $u=2,$ we have $B_{1}\leftarrow B_{2},$ so, no problem.
Further, we consider that $3\leq u\leq w.$ Let $i\in B_{u}.$ Let $k\in
\bigcup\limits_{t=1}^{u-1}B_{t}$ such that $P_{ik}>0.$ If $k\in B_{t}$ for
some $t\in \left\langle u-2\right\rangle ,$ then, using the definition of
sets $B_{1},B_{2},...,B_{w},$ we have $i\in B_{t+1}.$ Contradiction (because 
$t+1\leq u-1<u$ and $B_{u}\cap B_{t+1}=\emptyset $). So, $k\in B_{u-1},$ and
we have $B_{u-1}\leftarrow B_{u}.$ Using Theorem 2.8 for 
\[
B_{1}=\left\{ j\right\} \leftarrow B_{2}\leftarrow ...\leftarrow
B_{w-1}\leftarrow U_{0}=\left\{ i_{v}\right\} \leftarrow ...\leftarrow
U_{\left( n-b+1\right) s}=\left\langle n\right\rangle , 
\]
\noindent we have 
\[
\left( P^{w-1+\left( n-b+1\right) s}\right) ^{\left\{ j\right\} }>0, 
\]
\noindent where 
\[
w=\left\lfloor \frac{n-s-2}{k}\right\rfloor +3\text{ and }b=\max \left(
2,k\right) . 
\]

From Cases 1 and 2, since $\max \left( \left( n-b+1\right) s,w-1+\left(
n-b+1\right) s\right) =w-1+\left( n-b+1\right) s,$ by Theorem 4.1(i) we have 
\[
\left( P^{g}\right) ^{\left\{ j\right\} }>0, 
\]
\noindent where 
\[
g=w-1+\left( n-b+1\right) s=\left\lfloor \frac{n-s-2}{k}\right\rfloor
+2+s\left[ n-\max \left( 2,k\right) +1\right] \text{. }\square 
\]

\smallskip

Now, we give Wielandt Theorem and a generalization of it --- these are
corollaries of Theorems 4.22 and 4.23, respectively.

\smallskip

\textbf{THEOREM 4.24.} \textit{Let }$P\in N_{n},$\textit{\ }$n\geq 1$\textit{%
. Then }$P$\textit{\ is primitive if and only if} 
\[
P^{n^{2}-2n+2}>0. 
\]
\noindent (For this part --- this is Wielandt Theorem ---, see, \textit{e.g.}%
, [9, p. 520].)

\textit{More generally, if }$P\in \overline{G}_{1,1}^{+}\cup \overline{G}%
_{n,k}^{+},$\textit{\ }$n\geq 2$\textit{, }$k\in \left\langle
n-1\right\rangle ,$\textit{\ then }$P$\textit{\ is primitive if and only if} 
\[
P^{h}>0, 
\]

\noindent \textit{where} 
\[
h=\left\{ 
\begin{array}{ll}
1\medskip & \hspace{0.2cm}\text{\textit{if} }P\in \overline{G}_{1,1}^{+}, \\ 
n^{2}-2n+2\medskip & \hspace{0.2cm}\text{\textit{if} }P\in \overline{G}%
_{n,1}^{+}, \\ 
\left\lfloor \frac{n-m-3}{k}\right\rfloor +2+\left( m+1\right) \left(
n-k+1\right) & \hspace{0.2cm}\text{\textit{if} }P\in \overline{G}_{n,k}^{+},%
\text{ }k\neq 1,
\end{array}
\right. 
\]
\[
m=\left\lfloor \frac{n-2}{k}\right\rfloor +1\hspace{0.2cm}\text{\textit{if} }%
P\in \overline{G}_{n,k}^{+}. 
\]

\smallskip

\textit{Proof}. We consider the more general case. If $P\in \overline{G}%
_{1,1}^{+},$ then $P\in N_{1},$ $P>0.$ So, no problem. Further, we consider
that $P\in \overline{G}_{n,k}^{+},$ $n\geq 2,$ $k\in \left\langle
n-1\right\rangle .$

\smallskip

$``\Longrightarrow "$ Since $P\in \overline{G}_{n,k}^{+},$ we have $s\in
\left\langle m+1\right\rangle $, $s$ is that from Theorem 4.23 --- we prove
this statement. $P$ being irreducible, any $j\in \left\langle n\right\rangle 
$ belongs to a cycle. Fix $j\in \left\langle n\right\rangle .$ When $%
P_{jj}>0,$ we have $s=1$, and, therefore, $s\in \left\langle
m+1\right\rangle .$ When $P_{jj}=0,$ set 
\[
C_{0}=\left\{ j\right\} , 
\]
\[
C_{v}=\left\{ i\left| \text{ }i\in \left\langle n\right\rangle
-\bigcup\limits_{t=0}^{v-1}C_{t}\text{ and }\exists k\in
\bigcup\limits_{t=0}^{v-1}C_{t}\text{ such that }P_{ik}>0\right. \right\} , 
\]
\noindent $\forall v\geq 1$ with $\bigcup\limits_{t=0}^{v-1}C_{t}\subset
\left\langle n\right\rangle .$ $C_{0},$ $C_{1},$ $...,$ $C_{u}$ are disjoint
when $\bigcup\limits_{t=0}^{u-1}C_{t}\subset \left\langle n\right\rangle $
and $\bigcup\limits_{t=0}^{u}C_{t}=\left\langle n\right\rangle $ --- this
thing happens for some $u\geq 1.$ It follows that 
\[
\left| C_{0}\right| +\left| C_{1}\right| +...+\left| C_{u}\right| =n. 
\]
\noindent Since $P\in \overline{G}_{n,k}^{+},$ we have 
\[
\left| C_{1}\right| \geq k,\text{ }\left| C_{2}\right| \geq k,\text{ }...,%
\text{ }\left| C_{u-1}\right| \geq k,\text{ }\left| C_{u}\right| \geq a, 
\]
\noindent where $a=n-\left| \bigcup\limits_{t=0}^{u-1}C_{t}\right| $ ---
obviously, $1\leq a\leq k.$ So, 
\[
n\geq 1+\left( u-1\right) k+a 
\]
\noindent Further, we proceed as in the proof of Theorem 4.10, and obtain
that $u\leq m.$ Further, we have --- see the proof of Theorem 4.23, Case 2,
for a similar situation --- 
\[
C_{0}=\left\{ j\right\} \leftarrow C_{1}\leftarrow ...\leftarrow C_{u}. 
\]
\noindent Further, since $P$ is irreducible, $\exists q\in \left\langle
u\right\rangle $ ($q$ is unique or not; $q=u$ or $q\neq u$) such that 
\[
C_{q}\leftarrow \left\{ j\right\} . 
\]
\noindent So, $1\leq s\leq m+1.$ Finally, both when $P_{jj}>0$ and when $%
P_{jj}=0,$ we have $s\in \left\langle m+1\right\rangle $ --- the statement
was proved.

Fix $n$ and $k;$ $n\geq 2$\textit{, }$k\in \left\langle n-1\right\rangle .$
In this case, $g$ from Theorem 4.23 depends only on $s,$ and we consider the
function 
\[
g\left( s\right) =g=\left\lfloor \frac{n-s-2}{k}\right\rfloor +2+s\left[
n-\max \left( 2,k\right) +1\right] ,\text{ }\forall s\in \left\langle
m+1\right\rangle . 
\]

\smallskip

\noindent We show that $g$ is increasing --- moreover, it is strictly
increasing if $n\geq 3$ ($m\geq 2\Longrightarrow n\geq 3;$ $n\geq 3$ $\not%
{\Longrightarrow}$ $m\geq 2$). Fix $s\in \left\langle m\right\rangle .$

\smallskip

\textit{Case} 1. $k=1.$ We have 
\[
g\left( s+1\right) -g\left( s\right) =n-2\text{ }\left\{ 
\begin{array}{ll}
=0\smallskip & \text{if }n=2, \\ 
>0 & \text{if }n\geq 3.
\end{array}
\right. 
\]

\smallskip

\textit{Case} 2. $2\leq k\leq n-1.$ We have 
\[
g\left( s+1\right) -g\left( s\right) =n-k+1+\left\lfloor \frac{n-\left(
s+1\right) -2}{k}\right\rfloor -\left\lfloor \frac{n-s-2}{k}\right\rfloor . 
\]
\noindent Further, using the fact that $x-1<\left\lfloor x\right\rfloor \leq
x,$ $\forall x\in \Bbb{R}$, we have 
\[
g\left( s+1\right) -g\left( s\right) >n-k+1+\frac{n-\left( s+1\right) -2}{k}%
-1-\frac{n-s-2}{k}= 
\]
\[
=n-k-\frac{1}{k}>n-k-1\geq 0 
\]
\noindent (because $2\leq k\leq n-1$).

From Cases 1 and 2, a maximum value of $g$ is $g\left( m+1\right) ;$%
\[
g\left( m+1\right) =\left\lfloor \frac{n-m-3}{k}\right\rfloor +2+\left(
m+1\right) \left[ n-\max \left( 2,k\right) +1\right] . 
\]
\noindent If $s=m+1,$ no problem --- we have $P^{g\left( m+1\right) }>0.$ If 
$1\leq s\leq m,$ by Theorems 4.1(v) and 4.23 we have 
\[
P^{g\left( m+1\right) }=P^{g\left( m+1\right) -g\left( s\right) }P^{g\left(
s\right) }>0 
\]
\noindent ($g\left( s\right) =g,$ $g$ from Theorem 4.23, $P^{g\left(
s\right) }>0,$ ...).

In general, we cannot have $s\in \left\langle m\right\rangle $ --- see/study
for this fact the case when 
\[
P\in N_{3},\text{ }P=\left( 
\begin{array}{ccc}
0 & * & * \\ 
\ast & 0 & * \\ 
\ast & * & 0
\end{array}
\right) , 
\]
\noindent where ``$*$'' stands for a positive entry. But when $k=1,$ we have 
$s\in \left\langle m\right\rangle ,$ $m=n-1$ (if $s=n,$ we obtain that $P$
is periodic with period $n$ --- contradiction, see Theorem 4.21), and can
obtain a better result. For $k=1,$ we consider the function $g\!\mid
_{\left\langle m\right\rangle },$ the restriction of $g$ to $\left\langle
m\right\rangle .$ This function is also increasing, so, a maximum value of
it is $g\left( m\right) .$ Since $m=n-1,$ we obtain 
\[
g\left( n-1\right) =n^{2}-2n+2. 
\]
\noindent If $s=n-1$, no problem --- we have $P^{g\left( n-1\right) }>0.$ If 
$1\leq s\leq n-2,$ by Theorems 4.1(v) and 4.23 we have 
\[
P^{g\left( n-1\right) }=P^{g\left( n-1\right) -g\left( s\right) }P^{g\left(
s\right) }>0. 
\]

\smallskip

``$\Longleftarrow $'' Theorem 4.6. $\square $

\smallskip

The first upper bound, $n^{2}-2n+2,$ from Theorem 4.24 for the index of
primitivity is sharp --- for this, see the next example.

\smallskip

\textit{Example} 4.25. Consider the Wielandt matrix, 
\[
P\in N_{n},\text{ }P=\left( 
\begin{array}{cccccc}
0 & 1 & 0 & 0 & \cdots & 0 \\ 
0 & 0 & 1 & 0 & \cdots & 0 \\ 
\vdots & \vdots & \vdots & \ddots & \cdots & \vdots \\ 
0 & 0 & 0 & \cdots & \ddots & 0 \\ 
0 & 0 & 0 & 0 & \cdots & 1 \\ 
1 & 1 & 0 & 0 & \cdots & 0
\end{array}
\right) , 
\]

\noindent see, \textit{e.g.}, [9, Problem 4, p. 522]. We compute $\gamma
\left( P\right) $ using the $G^{+}$ method, more precisely, Theorem 2.8.
Since $\gamma \left( P\right) =\gamma \left( ^{\text{t}}P\right) ,$ we
compute $\gamma \left( ^{\text{t}}P\right) $ --- we have a reason, see
below. Set $Q=^{\text{t}}\!\!P.$ We consider that $Q=\left( Q_{ij}\right)
_{i,j\in \left\langle \left\langle n-1\right\rangle \right\rangle }$ (not $%
Q=\left( Q_{ij}\right) _{i,j\in \left\langle n\right\rangle }$) --- we made
this thing to use the addition modulo $n$. Set 
\[
U_{0}=\left\{ 0\right\} 
\]
\noindent and 
\[
U_{t}=\left\{ i\left| \text{ }i\in \left\langle \left\langle
n-1\right\rangle \right\rangle \text{ and }\exists j\in U_{t-1}\text{ such
that }Q_{ij}>0\right. \right\} ,\text{ }\forall t\geq 1. 
\]
\noindent We have 
\[
Q\in N_{n},\text{ }Q=\left( 
\begin{array}{cccccc}
0 & 0 & 0 & \cdots & 0 & 1 \\ 
1 & 0 & 0 & \cdots & 0 & 1 \\ 
0 & 1 & 0 & \cdots & 0 & 0 \\ 
\vdots & \vdots & \ddots & \cdots & \vdots & \vdots \\ 
0 & 0 & \cdots & \ddots & 0 & 0 \\ 
0 & 0 & 0 & \cdots & 1 & 0
\end{array}
\right) , 
\]
\noindent and 
\[
U_{0}=\left\{ 0\right\} \leftarrow U_{1}=\left\{ 1\right\} \leftarrow \ldots
\leftarrow U_{n-1}=\left\{ n-1\right\} \leftarrow ... 
\]

\noindent (here is the reason for the utilization of $^{\text{t}}P$--- to
have $\left\{ 0\right\} \leftarrow \left\{ 1\right\} \leftarrow \ldots
\leftarrow \left\{ n-1\right\} $) 
\[
\leftarrow U_{n}=U_{n+0\cdot \left( n-1\right) }=\left\{ 0,1\right\}
=\left\{ n\func{mod}n,\left( n+1\right) \func{mod}n\right\} \leftarrow 
\]
\[
\leftarrow U_{n+1}=\left\{ 1,2\right\} =\left\{ \left( n+1\right) \func{mod}%
n,\left( n+2\right) \func{mod}n\right\} \leftarrow ... 
\]

\[
...\leftarrow U_{n+\left( n-2\right) }=\left\{ n-2,n-1\right\} = 
\]

\[
=\left\{ \left( n+\left( n-2\right) \right) \func{mod}n,\left( n+\left(
n-1\right) \right) \func{mod}n\right\} \leftarrow 
\]
\[
\leftarrow U_{n+\left( n-1\right) }=U_{n+1\cdot \left( n-1\right) }=\left\{
n-1,0,1\right\} = 
\]
\[
=\left\{ \left( n+\left( n-1\right) \right) \func{mod}n,\left( n+n\right) 
\func{mod}n,\left( n+\left( n+1\right) \right) \func{mod}n\right\}
\leftarrow ... 
\]

\[
...\leftarrow U_{n+2\cdot \left( n-1\right) }\leftarrow ...\text{ }%
...\leftarrow U_{n+\left( n-2\right) \left( n-1\right) }=\left\langle
\left\langle n-1\right\rangle \right\rangle 
\]

\noindent --- $U_{n+0\cdot \left( n-1\right) }$ has $2$ elements, $%
U_{n+1\cdot \left( n-1\right) }$ has $3$ elements, and proceeding in this
way, we obtain that $U_{n+2\left( n-1\right) }$ has $4$ elements, ..., $%
U_{n+\left( n-2\right) \left( n-1\right) }$ has $n$ elements ($\left(
n-2\right) +2=n$), so, $U_{n+\left( n-2\right) \left( n-1\right)
}=\left\langle \left\langle n-1\right\rangle \right\rangle .$ Obviously,

\[
\left| U_{0}\right| \leq \left| U_{1}\right| \leq \ldots \leq \left|
U_{n+\left( n-2\right) \left( n-1\right) }\right| 
\]
\noindent and 
\[
\left| U_{n-1}\right| =1<\left| U_{n}\right| =\left| U_{n+0\cdot \left(
n-1\right) }\right| =2\leq ...\leq \left| U_{n+\left( n-1\right) -1}\right|
=2< 
\]
\[
<\left| U_{n+\left( n-1\right) }\right| =\left| U_{n+1\cdot \left(
n-1\right) }\right| =3\leq ...\text{ }... 
\]

\[
...\text{ }...\leq \left| U_{n+\left( n-2\right) \left( n-1\right)
-1}\right| =n-1<\left| U_{n+\left( n-2\right) \left( n-1\right) }\right| =n. 
\]

\noindent Since the number of sets $U_{t},$ $t\in \left\langle \left\langle
n+\left( n-2\right) \left( n-1\right) \right\rangle \right\rangle ,$ is 
\[
n+\left( n-2\right) \left( n-1\right) +1,\text{\textit{i.e.}},\text{ is }%
n^{2}-2n+3, 
\]
\noindent by Theorem 2.8 (see also Remark 2.9) we have 
\[
\left( Q^{n^{2}-2n+2}\right) ^{\left\{ 0\right\} }>0, 
\]
\[
\left( Q^{n^{2}-2n+2-1}\right) ^{\left\{ 1\right\} }=\left(
Q^{n^{2}-2n+1}\right) ^{\left\{ 1\right\} }>0,... 
\]
\[
...,\left( Q^{n^{2}-2n+2-\left( n-1\right) }\right) ^{\left\{ n-1\right\}
}=\left( Q^{n^{2}-3\left( n-1\right) }\right) ^{\left\{ n-1\right\} }>0. 
\]

\noindent We cannot have 
\[
\left( Q^{k}\right) ^{\left\{ 0\right\} }>0 
\]
\noindent for some $k\in \left\langle n^{2}-2n+1\right\rangle $ because $%
U_{t}\subset \left\langle \left\langle n-1\right\rangle \right\rangle ,$ $%
\forall t\in \left\langle \left\langle n^{2}-2n+2\right\rangle \right\rangle
.$ It follows that 
\[
Q^{n^{2}-2n+1}\ngtr 0\text{ and }Q^{n^{2}-2n+2}>0. 
\]
\noindent So, we have 
\[
\gamma \left( Q\right) =n^{2}-2n+2. 
\]

\smallskip

\textit{Problem} 4.26. (A good problem for the reader.) Is the second upper
bound, $h,$ from Theorem 4.24 for the index of primitivity sharp?

\smallskip

\textit{Remark} 4.27. (See, \textit{e.g.}, [9, p. 521].) To verify that a
matrix, $P,$ $P\in N_{n},$ is irreducible using Theorem 4.2, one way, a good
one, is, when $P\ngtr 0,$ to compute $\left( I+P\right) ^{2},$ $\left(
I+P\right) ^{4},$ $\left( I+P\right) ^{8},$ $...,$ $\left( I+P\right)
^{2^{t}},$ $t$ being the smallest integer for which $2^{t}\geq n-1.$ We can
proceed similarly, also when $P\ngtr 0,$ to verify that a matrix, $P,$ $P\in
N_{n},$ is primitive using Theorem 4.24. For more information, see, \textit{%
e.g.}, [9, p. 521]; see also Remark 4.4(b).

\bigskip

\begin{center}
\textbf{5. APPLICATIONS, II}
\end{center}

\bigskip

In this section, we give other applications of the $G^{+}$ method --- these
are also applications of Theorems 2.6 and 2.8. The results are also old and
new, and are for the irreducible matrices, for the reducible ones, for the
fully indecomposable ones, for the scrambling ones, and for the Sarymsakov
ones, in some cases the matrices being, moreover, $g_{k}^{+}$-matrices (not
only irreducible). Finally, we mention an application of Theorem 2.17.

\smallskip

For the first part of the next result, see, \textit{e.g.}, [14, p. 6] (in
this book, this part is considered as being important --- see, below, Remark
5.2).

\smallskip

\textbf{THEOREM 5.1.} \textit{Let }$y\in \Bbb{R}^{n},$\textit{\ }$n\geq 2,$%
\textit{\ be a nonnegative column vector with exactly }$h$\textit{\ positive
coordinates, }$h\in \left\langle n-1\right\rangle .$\textit{\ If }$P\in
N_{n},$\textit{\ }$n\geq 2$\textit{, is an irreducible matrix, then }$\left(
I+P\right) y$\textit{\ has more than }$h$\textit{\ positive coordinates.
More generally, if }$P\in \overline{G}_{n,k}^{+},$\textit{\ }$n\geq 2,$%
\textit{\ }$k\in \left\langle n-1\right\rangle ,$\textit{\ then }$\left(
I+P\right) y$\textit{\ has more than }$h+\min \left( k,n-h\right) -1$\textit{%
\ positive coordinates.}

\smallskip

\textit{Proof}. We prove the more general result. Let $Y\in N_{n,1}$
(therefore, $Y$ is a nonnegative column matrix), $Y_{i1}=y_{i},$ $\forall
i\in \left\langle n\right\rangle .$ Due this definition, we can work with $Y$
instead of $y$ --- we do this thing. Suppose that $%
Y_{i_{1}1},Y_{i_{2}1},...,Y_{i_{h}1}>0,$ where $i_{1},i_{2},...,i_{h}\in
\left\langle n\right\rangle ,$ $i_{u}\neq i_{v},$ $\forall u,v\in
\left\langle h\right\rangle ,$ $u\neq v.$ Set $V=\left\{
i_{1},i_{2},...,i_{h}\right\} .$ It follows that $Y\in \overline{G}%
_{V,\left\{ 1\right\} }^{+}.$

Now, we consider $I+P.$ $I+P\in \overline{G}_{n,k}^{+}$ because $P\in 
\overline{G}_{n,k}^{+}.$ Further, we have, by Theorem 3.2, $\emptyset \neq
D_{V}=D_{V}\left( I+P\right) \subseteq V^{c}$ and

\[
\left| D_{V}\right| \geq \min \left( k,\left| V^{c}\right| \right) =\min
\left( k,n-h\right) . 
\]

\noindent By Definition 3.1 (taking $E=D_{V}$), $I+P\in \overline{G}%
_{D_{V},V}^{+}.$ On the other hand, $I+P\in \overline{G}_{V,V}^{+}.$ So, $%
I+P\in \overline{G}_{V\cup D_{V},V}^{+}.$ By Theorem 2.6(i), $\left(
I+P\right) Y\in \overline{G}_{V\cup D_{V},\left\{ 1\right\} }^{+}.$ So, $%
\left( I+P\right) Y$ --- this is a nonnegative column matrix --- has $\left|
V\cup D_{V}\right| $ positive entries. It follows that $\left( I+P\right) Y$
has at least 
\[
h+\min \left( k,n-h\right) 
\]
\noindent positive entries because 
\[
\left| V\cup D_{V}\right| =\left| V\right| +\left| D_{V}\right| \geq h+\min
\left( k,n-h\right) 
\]
\noindent (recall that $\emptyset \neq D_{V}\subseteq V^{c};$ so, $D_{V}\cap
V=\emptyset $). $\square $

\smallskip

\textit{Remark} 5.2. (a) The above result can be used to give another proof
of Theorem 4.2, ``$\Longrightarrow $'', see, \textit{e.g.}, [14, p. 6] ---
[14] contains other interesting applications of Theorem 5.1 (see Corollary
2.1, p. 6, Theorems 2.2 and 2.3, p. 7, Corollary 4.2, p. 16, ...).

\smallskip

(b) Theorem 5.1 and the hypotheses of Theorems 4.8 and 4.10 suggest the
following result: if $y\in \Bbb{R}^{n},$\ $n\geq 2,$\ is a column vector
with exactly $h$\ positive coordinates, $h\in \left\langle n-1\right\rangle
, $\ and if $P\in N_{n},$\ $n\geq 2$, is an irreducible matrix and all the
main diagonal entries of $P$ are positive, then $Py$\ has more than $h$\
positive coordinates. More generally, if ... --- the completion and proof(s)
are left to the reader.

\smallskip

\textit{Definition} 5.3. (See, \textit{e.g.}, [8, pp. 34$-$35] and, better,
[14, p. 82].) Let $P\in N_{n},$ $n\geq 2$. Let $s\in \left\langle
n-1\right\rangle .$ We say that $P$ is \textit{partly decomposable} if it
contains an $s\times \left( n-s\right) $ zero submatrix.

\smallskip

\textit{Definition} 5.4. (See, \textit{e.g.}, [8, p. 35] and [14, pp. 82$-$%
83].) Let $P\in N_{n},$ $n\geq 2$. We say that $P$ is \textit{fully
indecomposable} if it is not partly decomposable.

\smallskip

By definition, the $1\times 1$ zero matrix is partly decomposable while a
nonzero $1\times 1$ matrix is fully indecomposable (see, \textit{e.g.}, [14,
p. 83]).

\smallskip

\textit{Remark} 5.5. The partly decomposable matrices are either irreducible
or reducible; they are generalizations of reducible matrices (any reducible
matrix is partly decomposable). The fully indecomposable matrices are
irreducible (not reducible), and, as a result, they are row-allowable.

\smallskip

\textit{Remark} 5.6. (See, \textit{e.g.}, [8, p. 35].) Let $P\in N_{n},$ $%
n\geq 2$. Then $P$ is fully indecomposable if and only if whenever it
contains a $p\times q$ zero submatrix, then $p+q\leq n-1.$

\smallskip

\textbf{THEOREM 5.7.} (See, \textit{e.g.}, [8, p. 37].)\textit{\ Let }$%
P_{1},P_{2},...,P_{n-1}\in N_{n},$\textit{\ }$n\geq 2$\textit{, be fully
indecomposable matrices. Then} 
\[
P_{1}P_{2}...P_{n-1}>0. 
\]

\smallskip

\textit{Proof}. The proof is less formal than that of Theorem 4.2, ``$%
\Longrightarrow $'', or that of Theorem 4.10. We show that 
\[
\left( P_{1}P_{2}...P_{n-1}\right) ^{\left\{ j\right\} }>0,\text{ }\forall
j\in \left\langle n\right\rangle . 
\]

Let $j\in \left\langle n\right\rangle .$

\textit{Case} 1. $\left( P_{n-1}\right) ^{\left\{ j\right\} }>0.$ When $n=2,$
we have only one matrix, $P_{1}$ --- nothing to prove. ($n=2$ and $P_{1}$ is
fully indecomposable $\Longrightarrow $ $P_{1}>0.$) When $n\geq 3,$ by
Theorem 4.1(i) (the fully indecomposable matrices are, by Remark 5.5 or
Remark 5.6, row-allowable), 
\[
\left( P_{1}P_{2}...P_{n-1}\right) ^{\left\{ j\right\} }=\left(
P_{1}P_{2}...P_{n-2}\right) \left( P_{n-1}\right) ^{\left\{ j\right\} }>0. 
\]

\textit{Case} 2. $\left( P_{n-1}\right) ^{\left\{ j\right\} }\ngtr 0.$ This
case holds when $n\geq 3$ (see above for $n=2$). Set 
\[
B_{1}=\left\{ j\right\} \text{ and }B_{2}=\left\{ i\left| \text{ }i\in
\left\langle n\right\rangle \text{ and }\left( P_{n-1}\right) _{ij}>0\right.
\right\} . 
\]
\noindent $B_{2}\neq \left\langle n\right\rangle $ because $\left(
P_{n-1}\right) ^{\left\{ j\right\} }\ngtr 0.$ It follows that $\left(
P_{n-1}\right) _{B_{2}^{c}}^{B_{1}}$ is a $\left| B_{2}^{c}\right| \times
\left| B_{1}\right| $ zero submatrix. By Remark 5.6 we have 
\[
\left| B_{2}^{c}\right| +\left| B_{1}\right| \leq n-1. 
\]
\noindent Further, we have 
\[
n-\left| B_{2}\right| +\left| B_{1}\right| \leq n-1, 
\]

\noindent and, therefore,

\[
\left| B_{2}\right| \geq \left| B_{1}\right| +1. 
\]

\noindent So, 
\[
\left| B_{2}\right| >\left| B_{1}\right| . 
\]

\noindent Moreover, $n>\left| B_{2}\right| >\left| B_{1}\right| $ (no
contradiction, because $n>2$). Set 
\[
B_{3}=\left\{ i\left| \text{ }i\in \left\langle n\right\rangle \text{ and }%
\left( P_{n-2}\right) _{ik}>0\text{ for some }k\in B_{2}\right. \right\} . 
\]

\noindent If $B_{3}=\left\langle n\right\rangle ,$ using Theorem 2.8 for $%
\left\{ j\right\} =B_{1}\leftarrow B_{2}\leftarrow B_{3}=\left\langle
n\right\rangle ,$ we have 
\[
\left( P_{n-2}P_{n-1}\right) ^{\left\{ j\right\} }>0. 
\]

\noindent If $n=3,$ no problem. If $n>3,$ by Theorem 4.1(i) we have 
\[
\left( P_{1}P_{2}...P_{n-1}\right) ^{\left\{ j\right\} }=\left(
P_{1}P_{2}...P_{n-3}\right) \left( P_{n-2}P_{n-1}\right) ^{\left\{ j\right\}
}>0 
\]

\noindent (by Remark 5.5, $P_{1},$ $P_{2},$ $...,$ $P_{n-3}$ are
row-allowable; $P_{1}P_{2}...P_{n-3}$ is row-allowable...). If $B_{3}\subset
\left\langle n\right\rangle ,$ proceeding as above, $\left( P_{n-2}\right)
_{B_{3}^{c}}^{B_{2}}$ is a $\left| B_{3}^{c}\right| \times \left|
B_{2}\right| $ zero submatrix... We proceed in this way until we obtain a
set $B_{v}$ with $B_{v}=\left\langle n\right\rangle ,$ where $v\geq 4.$ We
have $v\leq n$ because $1=\left| B_{1}\right| <\left| B_{2}\right|
<...<\left| B_{v}\right| =n.$ Using Theorem 2.8 for $\left\{ j\right\}
=B_{1}\leftarrow B_{2}\leftarrow ...\leftarrow B_{v}=\left\langle
n\right\rangle ,$ we have 
\[
\left( P_{n-v+1}P_{n-v+2}...P_{n-1}\right) ^{\left\{ j\right\} }>0. 
\]

\noindent If $v=n,$ no problem. If $v<n,$ by Theorem 4.1(i) we have 
\[
\left( P_{1}P_{2}...P_{n-1}\right) ^{\left\{ j\right\} }=\left(
P_{1}P_{2}...P_{n-v}\right) \left( P_{n-v+1}P_{n-v+2}...P_{n-1}\right)
^{\left\{ j\right\} }>0.\text{ }\square 
\]

\smallskip

The next result is a generalization of Theorem 5.7.

\smallskip

\textbf{THEOREM 5.8.} \textit{Let }$P_{1},P_{2},...,P_{m}\in \overline{G}%
_{n,k}^{+},$\textit{\ where} 
\[
m=\left\{ 
\begin{array}{ll}
n-1\smallskip & \text{if }k=1, \\ 
n-k+1 & \text{if }k\geq 2,
\end{array}
\right. 
\]
\noindent $n\geq 2$\textit{, and }$k\in \left\langle n-1\right\rangle .$%
\textit{\ If }$P_{1},P_{2},...,P_{m}$\textit{\ are fully indecomposable, then%
} 
\[
P_{1}P_{2}...P_{m}>0. 
\]

\smallskip

\textit{Proof}. For $k=1,$ see Theorem 5.7. Further, we suppose that $k\geq
2.$ The proof is similar to that of Theorem 5.7. Proceeding as there, for $%
j\in \left\langle n\right\rangle $ fixed, we have 2 cases, when $\left(
P_{m}\right) ^{\left\{ j\right\} }>0$ and when $\left( P_{m}\right)
^{\left\{ j\right\} }\ngtr 0.$ The first case is similar to Case 1 from the
proof of Theorem 5.7. The second one is similar to Case 2 from the proof of
Theorem 5.7 --- we construct/have the sets $B_{1},B_{2},...,B_{v},$ $%
B_{1}=\left\{ j\right\} \subset B_{2}\subset ...\subset B_{v}=\left\langle
n\right\rangle ,$ so, $1=\left| B_{1}\right| <\left| B_{2}\right|
<...<\left| B_{v}\right| =n.$ We have only a little difference: in the proof
of Theorem 5.7, Case 2, 
\[
\left| B_{2}\right| \geq \left| B_{1}\right| +1,\text{ }\left| B_{3}\right|
\geq \left| B_{2}\right| +1,\text{ }...,\text{ }\left| B_{v}\right| \geq
\left| B_{v-1}\right| +1, 
\]

\noindent while here, since $P_{1},P_{2},...,P_{m}\in \overline{G}%
_{n,k}^{+}, $ we have 
\[
\left| B_{2}\right| \geq \max \left( k,\left| B_{1}\right| +1\right) =\max
\left( k,2\right) =k, 
\]
\[
\left| B_{3}\right| \geq \max \left( k,\left| B_{2}\right| +1\right) \geq
\max \left( k,k+1\right) =k+1, 
\]
\[
\vdots 
\]
\[
\left| B_{n-k+2}\right| \geq k+\left( n-k\right) =n 
\]

\noindent --- therefore, $v\leq n-k+2,$ and the number of matrices we need
(the worst case) is $n-k+1,$ \textit{i.e.}, is $m$ (we justified the
definition of $m$ for $k\geq 2$). $\square $

\smallskip

By Theorem 5.8, if $P\in \overline{G}_{n,k}^{+}$ and is fully
indecomposable, then --- it is interesting this thing --- 
\[
\gamma \left( P\right) \leq \left\{ 
\begin{array}{ll}
n-1\smallskip & \text{if }k=1, \\ 
n-k+1 & \text{if }k\geq 2.
\end{array}
\right. 
\]

\smallskip

The nonnegative matrices which have at least one positive column are of
interest, \textit{e.g.}, in the finite Markov chain theory, and, in this
case, the matrices are stochastic. If $P\in S_{r}$ ($P$ can be considered as
being the transition matrix of a Markov chain) and has at least one positive
column, then $\mu \left( P\right) >0,$ where 
\[
\mu \left( P\right) =\max_{j\in \left\langle r\right\rangle }\left(
\min_{i\in \left\langle r\right\rangle }P_{ij}\right) , 
\]

\noindent $\mu \left( P\right) $ is the Markov ergodicity coefficient of $P$
(see, \textit{e.g.}, [10, pp. 56$-$57]), and, as a result --- this is
important ---, $\lim\limits_{n\rightarrow \infty }P^{n}$ exists. \textit{E.g.%
}, 
\[
P=\left( 
\begin{array}{cc}
0 & 1 \\ 
\frac{1}{2} & \frac{1}{2}
\end{array}
\right) 
\]

\noindent has a positive column (moreover, it is a primitive matrix, and a
main diagonal entry of it is $0$), and 
\[
\lim\limits_{n\rightarrow \infty }P^{n}=\left( 
\begin{array}{cc}
\frac{1}{3}\smallskip & \frac{2}{3} \\ 
\frac{1}{3} & \frac{2}{3}
\end{array}
\right) . 
\]

\noindent Another example, when 
\[
P=\left( 
\begin{array}{cc}
1 & 0 \\ 
\frac{1}{2} & \frac{1}{2}
\end{array}
\right) 
\]
\noindent ($P$ has a positive column; moreover, it is reducible and all the
main diagonal entries of it are positive), is left to the reader. For more
information, see, \textit{e.g.}, [10, Section 1.11 and Chapter 4]; see, 
\textit{e.g.}, also [8, pp. 61$-$62].

\smallskip

\textit{Definition} 5.9. (See, \textit{e.g.}, [10, p. 57] and [21, p. 140].)
Let $P\in N_{m,n},$ $m,n\geq 1$. We say that $P$ is a \textit{Markov matrix}
if it has at least one positive column (equivalently, if $P$ is sum-positive
on $\left\langle m\right\rangle \times \left\{ j\right\} $ for some $j\in
\left\langle n\right\rangle $).

\smallskip

Theorem 2.8 is a general result for the products of nonnegative matrices,
products which are Markov --- it is also a general result to show that a
product of nonnegative matrices is positive. Below we give a few results for
the products of nonnegative matrices, products which have a positive column
(therefore, are Markov or more).

\smallskip

\textbf{THEOREM 5.10.} (See, \textit{e.g.}, [8, p. 62] --- our result is
better.) \textit{Let }$P\in N_{n},$\textit{\ }$n\geq 2$\textit{, be in the
partitioned form} 
\[
P=\left( 
\begin{array}{cc}
Q & 0 \\ 
R & T
\end{array}
\right) , 
\]
\textit{\noindent where }$Q\in N_{m}$\textit{\ and is a primitive matrix }($%
1\leq m<n$). \textit{If for any }$i\in \left\{ m+1,m+2,...,n\right\} $ ($%
\left\{ m+1,m+2,...,n\right\} =\left\langle m\right\rangle ^{c}$),\textit{\
there exists a path from }$i$\textit{\ to some, say, }$j,$\textit{\ }$j\in
\left\langle m\right\rangle ,$ $j=j\left( i\right) $\textit{\ }(\textit{see
Definition }4.17), \textit{then} 
\[
P^{\gamma \left( Q\right) \left( n-m+1\right) } 
\]
\textit{\noindent has the first }$m$\textit{\ columns positive }(\textit{%
equivalently,} 
\[
\left( P^{\gamma \left( Q\right) \left( n-m+1\right) }\right) ^{\left\langle
m\right\rangle }>0\text{),} 
\]

\noindent \textit{and, more generally,} 
\[
\left( P^{t\left( n-m+1\right) }\right) ^{\left\langle m\right\rangle }>0%
\text{, }\forall t\geq \gamma \left( Q\right) 
\]

\noindent ($\gamma \left( Q\right) $ \textit{is the index of primitivity of} 
$Q$).

\smallskip

\textit{Proof}. Since $Q^{\gamma \left( Q\right) }>0,$ we have 
\[
\left( P^{\gamma \left( Q\right) }\right) _{\left\langle m\right\rangle
}^{\left\langle m\right\rangle }>0. 
\]

\noindent Let $Z=P^{\gamma \left( Q\right) }.$ $Z$ also has the above
property --- a path property: $\forall i\in \left\langle m\right\rangle
^{c}, $ there exists a path from $i$ to some, say, $k,$ $k\in \left\langle
m\right\rangle ,$ $k=k\left( i\right) $ (an elementary proof --- use
Definition 4.17 and the definition of $P$). We show that $\forall j\in
\left\langle m\right\rangle ,$ $\exists l=l\left( j\right) \in \left\langle
n-m+1\right\rangle $ such that 
\[
\left( Z^{l}\right) ^{\left\{ j\right\} }>0. 
\]

Let $j\in \left\langle m\right\rangle $. Set 
\[
D_{1}=\left\{ j\right\} . 
\]
\noindent If $Z^{\left\{ j\right\} }>0,$ no problem, we take $l=1.$ If $%
Z^{\left\{ j\right\} }\ngtr 0,$ we consider a new set, 
\[
D_{2}=D_{1}\cup \left\{ i\left| \text{ }i\in \left\langle n\right\rangle 
\text{ and }Z_{ij}>0\right. \right\} . 
\]
\noindent Since $Z_{ij}>0,$ $\forall i\in \left\langle m\right\rangle ,$ we
have $D_{2}\supseteq \left\langle m\right\rangle $ (not $\supset
\left\langle m\right\rangle ,$ because we can have $Z_{ij}=0,$ $\forall i\in
\left\langle m\right\rangle ^{c}$). Since $Z^{\left\{ j\right\} }\ngtr 0,$
we have $D_{2}\subset \left\langle n\right\rangle .$ So, $\left\langle
n\right\rangle \supset D_{2}\supseteq \left\langle m\right\rangle .$
Further, we set 
\[
D_{3}=D_{2}\cup \left\{ i\left| \text{ }i\in \left\langle n\right\rangle 
\text{ and }\exists k\in D_{2}\text{ such that }Z_{ik}>0\right. \right\} . 
\]
\noindent Obviously, $D_{3}\supset D_{2}$ because $\left\langle
n\right\rangle \supset D_{2}\supseteq \left\langle m\right\rangle $ and $Z$
has the above path property. If $D_{3}=\left\langle n\right\rangle ,$ using
Theorem 2.8, we have $\left( Z^{2}\right) ^{\left\{ j\right\} }>0,$ and we
take $l=2.$ If $D_{3}\subset \left\langle n\right\rangle ,$ we construct a
new set, 
\[
D_{4}=D_{3}\cup \left\{ i\left| \text{ }i\in \left\langle n\right\rangle 
\text{ and }\exists k\in D_{3}\text{ such that }Z_{ik}>0\right. \right\} . 
\]

\noindent Obviously, $D_{4}\supset D_{3}$ because $\left\langle
n\right\rangle \supset D_{3}\supset D_{2}\supseteq \left\langle
m\right\rangle $ and $Z$ has the above path property. If ... We proceed in
this way until we obtain a set $D_{u},$ $u=u\left( j\right) \geq 1,$ with $%
D_{u}=\left\langle n\right\rangle .$ Obviously, $1\leq u\leq 2+n-m$ because $%
\left\langle m\right\rangle \subseteq D_{2}\subset D_{3}\subset ...\subset
D_{u}=\left\langle n\right\rangle .$ In this case (when we have $%
D_{u}=\left\langle n\right\rangle $), using Theorem 2.8 for $D_{1}\leftarrow
D_{2}\leftarrow ...\leftarrow D_{u}$, we have $\left( Z^{u-1}\right)
^{\left\{ j\right\} }>0,$ and we take $l=u-1.$ It follows that $l\leq n-m+1.$

From the above results, by Theorem 4.1(i) ($P$ is row-allowable $%
\Longrightarrow $ $\forall t\geq 1,$ $P^{t}$ is row-allowable $%
\Longrightarrow $ $\forall t\geq 1,$ $Z^{t}$ is row-allowable) 
\[
\left( Z^{n-m+1}\right) ^{\left\langle m\right\rangle }>0\text{.} 
\]
\noindent So, 
\[
\left( P^{\gamma \left( Q\right) \left( n-m+1\right) }\right) ^{\left\langle
m\right\rangle }>0\text{,} 
\]

\noindent and, more generally, 
\[
\left( P^{t\left( n-m+1\right) }\right) ^{\left\langle m\right\rangle }>0%
\text{, }\forall t\geq \gamma \left( Q\right) .\text{ }\square 
\]

\smallskip

\textit{Definition} 5.11. (See, \textit{e.g.}, [8, p. 63], [10, p. 57], and
[21, p. 143].) Let $P\in N_{m,n},$ $m\geq 2,$ $n\geq 1$. We say that $P$ is
a \textit{scrambling matrix} if $\forall i,j\in \left\langle m\right\rangle
, $ $i\neq j,$ $\exists k\in \left\langle n\right\rangle $ such that $%
P_{ik},P_{jk}>0$ (equivalently, if $\forall i,j\in \left\langle
m\right\rangle ,$ $i\neq j,$ $P$ is sum-positive on $\left\{ i,j\right\}
\times \left\{ k\right\} $ for some $k\in \left\langle n\right\rangle $).

\smallskip

\textbf{THEOREM 5.12.} (See, e.g.\textit{, }[8, p. 63] --- our result is
better.) \textit{Let }$P\in N_{n},$\textit{\ }$n\geq 2$\textit{. If }$P$%
\textit{\ is a scrambling matrix, then} 
\[
P^{n-1} 
\]
\noindent \textit{has a positive column, i.e., is a Markov matrix.}

\smallskip

\textit{Proof}. \textit{Case} 1. $P$ has a positive column. (Obviously, in
this case, $P$ is a scrambling matrix.) It follows that $P^{n-1}$ has a
positive column too because the scrambling matrices are row-allowable.

\smallskip

\textit{Case} 2. $P$ has no positive columns. (Obviously, there exist
scrambling matrices which have no positive columns.) We do a head-to-tail
construction for Theorem 2.8 (to apply Theorem 2.8). Set 
\[
U_{1}=\left\langle n\right\rangle \text{ and }t_{1}=\left| U_{1}\right| =n. 
\]
\noindent Let $j_{1}^{\left( 1\right) },$ $j_{2}^{\left( 1\right) },$ $...,$ 
$j_{n-1}^{\left( 1\right) }\in \left\langle n\right\rangle $ such that 
\[
P_{1j_{1}^{\left( 1\right) }},\text{ }P_{nj_{1}^{\left( 1\right) }}>0,\text{ 
}P_{2j_{2}^{\left( 1\right) }},\text{ }P_{nj_{2}^{\left( 1\right) }}>0,\text{
}...,\text{ }P_{n-1\rightarrow j_{n-1}^{\left( 1\right) }},\text{ }%
P_{nj_{n-1}^{\left( 1\right) }}>0 
\]
\noindent --- this thing can be done because $P$ is a scrambling matrix. Set 
\[
U_{2}=\bigcup_{l=1}^{n-1}\left\{ j_{l}^{\left( 1\right) }\right\} \text{ and 
}t_{2}=\left| U_{2}\right| . 
\]
\noindent We have $t_{2}>1$ because $P$ has no positive columns. Obviously, $%
U_{2}\subset \left\langle n\right\rangle .$ So, $1<t_{2}<t_{1}.$ Consider
that 
\[
U_{2}=\left\{ i_{1}^{\left( 2\right) },i_{2}^{\left( 2\right)
},...,i_{t_{2}}^{\left( 2\right) }\right\} . 
\]
\noindent Let $j_{1}^{\left( 2\right) },$ $j_{2}^{\left( 2\right) },$ $...,$ 
$j_{t_{2}-1}^{\left( 2\right) }\in \left\langle n\right\rangle $ such that

\[
P_{i_{1}^{\left( 2\right) }j_{1}^{\left( 2\right) }},\text{ }%
P_{i_{t_{2}}^{\left( 2\right) }j_{1}^{\left( 2\right) }}>0,\text{ }%
P_{i_{2}^{\left( 2\right) }j_{2}^{\left( 2\right) }},\text{ }%
P_{i_{t_{2}}^{\left( 2\right) }j_{2}^{\left( 2\right) }}>0,\text{ }...,\text{
}P_{i_{t_{2}}^{\left( 2\right) }-1\rightarrow j_{t_{2}-1}^{\left( 2\right)
}},\text{ }P_{i_{t_{2}}^{\left( 2\right) },j_{t_{2}-1}^{\left( 2\right)
}}>0. 
\]

\noindent Set 
\[
U_{3}=\bigcup_{l=1}^{t_{2}-1}\left\{ j_{l}^{\left( 2\right) }\right\} \text{
and }t_{3}=\left| U_{3}\right| . 
\]
\noindent Obviously, we have $t_{2}>t_{3}\geq 1.$ Consider that 
\[
U_{3}=\left\{ i_{1}^{\left( 3\right) },i_{2}^{\left( 3\right)
},...,i_{t_{3}}^{\left( 3\right) }\right\} . 
\]

\noindent If $t_{3}=1,$ since $\left\langle n\right\rangle =U_{1}\rightarrow
U_{2}\rightarrow U_{3},$ by Theorem 2.8 we have 
\[
\left( P^{2}\right) ^{\left\{ i_{1}^{\left( 3\right) }\right\} }>0. 
\]

\noindent If $t_{3}>1,$ let $j_{1}^{\left( 3\right) },$ $j_{2}^{\left(
3\right) },$ $...,$ $j_{t_{3}-1}^{\left( 3\right) }\in \left\langle
n\right\rangle $ such that... We proceed in this way until we find a set $%
U_{v}$ with $\left| U_{v}\right| :=t_{v}=1$ --- this thing happens because $%
n=t_{1}>t_{2}>...>t_{v},$ and, obviously, $2\leq v\leq n.$ By Theorem 2.8,
since $\left\langle n\right\rangle =U_{1}\rightarrow U_{2}\rightarrow
...\rightarrow U_{v}$ and $\left| U_{v}\right| =1,$ we have 
\[
\left( P^{v-1}\right) ^{U_{v}}>0. 
\]

\noindent If $v=n,$ no problem. If $2\leq v<n$, by Theorem 4.1(iii) ($%
P^{n-1}=P^{v-1}P^{n-v},$ $P^{n-v}$ is row-allowable, ...), $\exists k\in
\left\langle n\right\rangle $ such that 
\[
\left( P^{n-1}\right) ^{\left\{ k\right\} }>0, 
\]

\noindent \textit{i.e.}, $P^{n-1}$ is a Markov matrix. $\square $

\smallskip

Theorem 5.12 can be generalized (see also [8, p. 64, Theorem 4.6] for a less
good result).

\smallskip

\textbf{THEOREM 5.13.} \textit{Let }$P_{1},P_{2},...,P_{n-1}\in N_{n},$%
\textit{\ }$n\geq 2$\textit{. If }$P_{1},P_{2},...,P_{n-1}$\textit{\ are
scrambling matrices, then} 
\[
P_{1}P_{2}...P_{n-1} 
\]
\noindent \textit{is a Markov matrix. More generally, if }$P_{1}\in
N_{n_{1},n_{2}},$ $P_{2}\in N_{n_{2},n_{3}},$ $...,$ $P_{t}\in
N_{n_{t},n_{t+1}},$\textit{\ }$n_{1},n_{2},...,n_{t}\geq 2,$ $n_{t+1}\geq 1,$
$t=n_{1}-1,$ and $P_{1},P_{2},...,P_{t}$\textit{\ are scrambling matrices,
then} 
\[
P_{1}P_{2}...P_{t} 
\]
\noindent \textit{is a Markov matrix --- we even have a better result,
namely,} 
\[
P_{1}P_{2}...P_{z} 
\]
\noindent \textit{is a Markov matrix, where }$z=\min\limits_{0\leq u\leq
t-1}\left( u+n_{u+1}-1\right) .$

\smallskip

\textit{Proof}. Similar to the proof of Theorem 5.12$.$ As to $z$, $z\leq t$
because $u+n_{u+1}-1$ is equal to $t$ if $u=0;$ we have $u$ matrices before
the matrix $P_{u+1},$ and when $u+n_{u+1}-1\leq t,$ the number of matrices $%
P_{u+1},$ $P_{u+2},$ $...,$ $P_{u+n_{u+1}-1}$ is $n_{u+1}-1,$ so, when $%
u+n_{u+1}-1\leq t,$ we have $u+n_{u+1}-1$ matrices; the cases when $%
u+n_{u+1}-1>t$ do not count because $z\leq t$. $\square $

\smallskip

Let $P\in N_{m,n},$ $m,n\geq 1$. Let $\emptyset \neq T\subseteq \left\langle
m\right\rangle .$ Set 
\[
F\left( T\right) =\left\{ j\left| \text{ }j\in \left\langle n\right\rangle 
\text{ and for which }\exists i\in T\text{ such that }P_{ij}>0\right.
\right\} . 
\]
\noindent We call $F\left( T\right) $ the \textit{set of consequent indices
of }(\textit{set of indices})\textit{\ }$T$ (see, \textit{e.g.}, [8, p. 64]
and [21, p. 146]).

\smallskip

\textit{Definition} 5.14. (See, \textit{e.g.}, [8, p. 65], [21, p. 146], and
[23].) Let $P\in N_{m,n},$ $m\geq 2,$ $n\geq 1$. We say that $P$ is a 
\textit{Sarymsakov matrix} if $\forall I,J,$ $\emptyset \neq I,J\subseteq
\left\langle m\right\rangle ,$ $I\cap J=\emptyset ,$ either

(1) $F\left( I\right) \cap F\left( J\right) \neq \emptyset $

\noindent or

(2) $F\left( I\right) \cap F\left( J\right) =\emptyset $ and $\left| F\left(
I\right) \cup F\left( J\right) \right| >\left| I\cup J\right| .$

\smallskip

\textit{Remark} 5.15. (a) If $P$ is a Sarymsakov matrix, then it is a
row-allowable matrix (see, \textit{e.g.}, [8, p. 65]).

(b) If $P$ is a scrambling matrix, then it is a Sarymsakov matrix (see, 
\textit{e.g.}, [21, p. 146]).

\smallskip

\textbf{THEOREM 5.16.} \textit{Let }$P\in N_{m,n}$\textit{\ and }$Q\in
N_{n,p},$\textit{\ }$m\geq 2,$\textit{\ }$n,p\geq 1$\textit{. If }$P$\textit{%
\ is sum-positive on }$\left\{ i,j\right\} \times \left\{ k\right\} $\textit{%
\ }($\left\{ i,j\right\} \subseteq \left\langle m\right\rangle ,$\textit{\ }$%
\left\{ k\right\} \subseteq \left\langle n\right\rangle $)\textit{\ and }$Q$%
\textit{\ is row-allowable or, more generally, }$Q_{\left\{ k\right\} }$%
\textit{\ is row-allowable, then }$PQ$\textit{\ is sum-positive on }$\left\{
i,j\right\} \times \left\{ l\right\} $\textit{\ for some }$l\in \left\langle
p\right\rangle $\textit{.}

\smallskip

\textit{Proof}. Since $P$ is sum-positive on $\left\{ i,j\right\} \times
\left\{ k\right\} $ and $Q$ is sum-positive on $\left\{ k\right\} \times
\left\{ l\right\} $ for some $l\in \left\langle p\right\rangle $, then, by
Theorem 2.6(i), $PQ$ is sum-positive on $\left\{ i,j\right\} \times \left\{
l\right\} $ for some $l\in \left\langle p\right\rangle .$ $\square $

\smallskip

\textbf{THEOREM 5.17.} (See, \textit{e.g.}, [8, p. 66] and [21, p. 146].)%
\textit{\ Let }$P_{1},$ $P_{2},$ $...,$ $P_{n-1}\in N_{n},$\textit{\ }$n\geq
2.$\textit{\ If }$P_{1},$ $P_{2},$ $...,$ $P_{n-1}$\textit{\ are Sarymsakov
matrices, then} 
\[
P_{1}P_{2}...P_{n-1} 
\]
\noindent \textit{is a scrambling matrix.}

\smallskip

\textit{Proof}. See, \textit{e.g.}, [8, p. 66] and [21, p. 146] --- Remark
5.15(a) and Theorem 5.16 (not using our terminology) are used. $\square $

\smallskip

\textit{Remark} 5.18. Theorem 5.17 can be generalized considering $P_{1}\in
N_{m_{1},m_{2}},$ $P_{2}\in N_{m_{2},m_{3}},$ $...,$ $P_{n-1}\in
N_{m_{n-1},m_{n}},$ $m_{1},$ $m_{2},$ $...,$ $m_{n-1}\geq 2,$ $m_{n}\geq 1,$ 
$m_{1},$ $m_{2},$ $...,$ $m_{n-1}\leq n$, and $P_{1},$ $P_{2},$ $...,$ $%
P_{n-1}$ be Sarymsakov matrices.

\smallskip

Recall that one of our aim is to obtain Markov matrices --- these are of
interest, \textit{e.g.}, in the finite Markov chain theory. It is also of
interest in the finite Markov chain theory to obtain scrambling matrices,
see, \textit{e.g.}, Theorem 5.17 for an example, because, if, \textit{e.g.},
a stochastic $r\times r$ matrix, $r\geq 2,$ say, $P,$ is scrambling, then $%
\alpha \left( P\right) >0,$ where 
\[
\alpha \left( P\right) =\min_{1\leq i,j\leq r}\sum_{k=1}^{r}\min \left(
P_{ik},P_{jk}\right) , 
\]
\noindent $\alpha \left( P\right) $ is the Dobrushin ergodicity coefficient
of $P$ (see, \textit{e.g.}, [10, pp. 56$-$57]), and, as a result --- this is
important, --- $\lim\limits_{n\rightarrow \infty }P^{n}$ exists.

\smallskip

\textbf{THEOREM 5.19.}\textit{\ Let }$P_{1},$ $P_{2},$ $...,$ $P_{\left(
n-1\right) ^{2}}\in N_{n},$\textit{\ }$n\geq 2.$\textit{\ If }$P_{1},$ $%
P_{2},$ $...,$ $P_{\left( n-1\right) ^{2}}$\textit{\ are Sarymsakov
matrices, then} 
\[
P_{1}P_{2}...P_{\left( n-1\right) ^{2}} 
\]
\noindent \textit{is a Markov matrix. In particular, if }$%
P_{1}=P_{2}=...=P_{\left( n-1\right) ^{2}}:=P,$ \textit{then} 
\[
P^{\left( n-1\right) ^{2}} 
\]
\noindent \textit{is a Markov matrix.}

\smallskip

\textit{Proof}. Theorems 5.13 and 5.17. $\square $

\smallskip

Recall that the $G^{+}$ method was suggested by the $G$ method and Theorem
2.2, both from [17]. Theorem 2.2 from [17] can be proved --- it is easy ---
using Theorem 2.17 from here. So, we have an application of Theorem 2.17, an
important one because Theorem 2.2 from [17] was used to show that the
transition matrix of our hybrid Metropolis-Hastings chain from [17] is
positive, see Theorem 2.3 in [17], see also Theorem 1.2 in [20] ---
interestingly, this positive matrix is a product of stochastic matrices, all
being reducible or, excepting the first one, the others are reducible. Since
our Gibbs sampler in a generalized sense (see, \textit{e.g.}, [20] for this
chain) is a special case, an important one, of our hybrid
Metropolis-Hastings chain, special cases for Theorem 2.3 from [17] can be
found where we applied our Gibbs sampler in a generalized sense, see, 
\textit{e.g.}, [18] and [19] --- in each application from there, the
transition matrix of our Gibbs sampler in a generalized sense is a product
of reducible stochastic matrices; moreover, this transition matrix is stable
(for stable matrices, see, \textit{e.g.}, [19]). Another special case for
Theorem 2.3 from [17] is at the cyclic Gibbs sampler in the finite case
because this chain is a special case of our Gibbs sampler in a generalized
sense (and, therefore, a special case of our hybrid Metropolis-Hastings
chain).

\smallskip

Other applications, important applications, of Theorem 2.17 will be found
--- we believe this --- in the future.

\smallskip

At present we have two $G$-type methods, the $G$ method and $G^{+}$ method,
... and the results obtained using them are impressive.

\bigskip

\begin{center}
\textbf{REFERENCES}
\end{center}

\bigskip \ 

[1] J. Bang-Jensen and G.Z. Gutin, \textit{Digraphs}:\textit{\ Theory,
Algorithms and Applications}, 2nd Edition. Springer, London, 2009.

[2] A. Benjamin, G. Chartrand, and P. Zhang, \textit{The Fascinating World
of Graph Theory}. Princeton University Press, Princeton, 2015.

[3] J.A. Bondy and U.S.R. Murty, \textit{Graph Theory with Applications}.
North-Holland, New York, 1979.

[4] J.A. Bondy and U.S.R. Murty, \textit{Graph Theory}. Springer, London,
2008.

[5] S. Brooks, A. Gelman, G.L. Jones, and X.-L. Meng (Eds.), \textit{%
Handbook of Markov Chain Monte Carlo}. Chapman \& Hall/CRC, Boca Raton, 2011.

[6] R. Diestel, \textit{Graph Theory}, 3rd Edition. Springer, Heidelberg,
2005.

[7] J. Ding and A. Zhou, \textit{Nonnegative Matrices, Positive Operators,
and Applications}. World Scientific, New Jersey, 2009.

[8] D.J. Hartfiel, \textit{Nonhomogeneous Matrix Products}. World
Scientific, Singapore, 2002.

[9] R.A. Horn and C.R. Johnson, \textit{Matrix Analysis}. Cambridge
University Press, Cambridge, 1985.

[10] M. Iosifescu, \textit{Finite Markov Processes and Their Applications}.
Wiley, Chichester \& Ed. Tehnic\u {a}, Bucharest, 1980; corrected
republication by Dover, Mineola, N.Y., 2007.

[11] G.F. Lawler and L.N. Coyle, \textit{Lectures on Contemporary Probability%
}. AMS, Providence, Rhode Island, 1999.

[12] D.A. Levin, Y. Peres, and E.L. Wilmer, \textit{Markov Chains and Mixing
Times}. AMS, Providence, Rhode Island, 2009.

[13] N. Madras and G. Slade, \textit{The Self-Avoiding Walk}.
Birkh\"{a}user, Boston, 1996.

[14] H. Minc, \textit{Nonnegative Matrices}. Wiley, New York, 1988.

[15] A. Paz, \textit{Introduction to Probabilistic Automata}. Academic
Press, New York, 1971.

[16] U. P\u {a}un, $G_{\Delta _{1},\Delta _{2}}$ \textit{in action}. Rev.
Roumaine Math. Pures Appl.\textbf{\ 55 }(2010), 387$-$406.

[17] U. P\u {a}un, \textit{A hybrid Metropolis-Hastings chain}. Rev.
Roumaine Math. Pures Appl.\textbf{\ 56 }(2011), 207$-$228.

[18] U. P\u {a}un, $G$\textit{\ method in action}: \textit{fast exact
sampling from set of permutations of order }$n$\textit{\ according to
Mallows model through Kendall metric.} Rev. Roumaine Math. Pures Appl.%
\textbf{\ 63 }(2018), 259$-$280.

[19] U. P\u {a}un, \textit{Ewens distribution on }$\Bbb{S}_{n}$ \textit{is a
wavy probability distribution with respect to }$n$\textit{\ partitions.} An.
Univ. Craiova Ser. Mat. Inform. \textbf{47} (2020), 1$-$24.

[20] U. P\u {a}un, $\Delta $-\textit{wavy probability distributions and
Potts model.} An. Univ. Craiova Ser. Mat. Inform. \textbf{49} (2022), 208$-$%
249.

[21] E. Seneta, \textit{Non-negative Matrices and Markov Chains}, 2nd
Edition. Springer-Verlag, Berlin, 1981; revised printing, 2006.

[22] W.T. Tutte, \textit{Graph Theory}. Cambridge University Press,
Cambridge, 2001.

[23] W. Xia, J. Liu, M. Cao, K.H. Johansson, and T. Ba\c{s}ar, \textit{%
Generalized Sarymsakov matrices}. IEEE Trans. Automat. Control \textbf{64}
(2019), 3085$-$3100.

\bigskip \ 
\[
\begin{array}{ccc}
\mathit{April}\text{ }\mathit{7,}\text{ }\mathit{2023} &  & \mathit{%
Romanian\ Academy} \\ 
&  & \mathit{Gheorghe\ Mihoc}-\mathit{Caius\ Iacob}\text{ }\mathit{Institute}
\\ 
&  & \mathit{of\ Mathematical\ Statistics}\text{ }\mathit{and\ Applied\
Mathematics} \\ 
&  & \mathit{Calea\ 13\ Septembrie\ nr.\ 13} \\ 
&  & \mathit{050711\ Bucharest\ 5,\ Romania} \\ 
&  & \mathit{upterra@gmail.com}
\end{array}
\]

\end{document}